\documentclass[12pt]{amsart}
\usepackage[margin=1.2in]{geometry}
\usepackage[T1]{fontenc}
\usepackage{amsmath,amssymb,amsthm,mathtools,mathrsfs}
\usepackage{microtype,booktabs,array,graphicx,tikz}
\usetikzlibrary{arrows.meta,calc}
\usepackage[colorlinks=true,linkcolor=blue,citecolor=blue,urlcolor=blue]{hyperref}
\numberwithin{equation}{section}
\newtheorem{theorem}[equation]{Theorem}
\newtheorem{proposition}[equation]{Proposition}
\newtheorem{lemma}[equation]{Lemma}
\newtheorem{corollary}[equation]{Corollary}
\theoremstyle{definition}
\newtheorem{definition}[equation]{Definition}
\newtheorem{example}[equation]{Example}
\theoremstyle{remark}
\newtheorem{remark}[equation]{Remark}
\newcommand{\R}{\mathbb R}
\newcommand{\FM}{\operatorname{FM}}
\newcommand{\End}{\operatorname{End}}
\newcommand{\Conv}{\operatorname{Conv}}
\newcommand{\sgn}{\operatorname{sgn}}
\newcommand{\gA}{\mathfrak g_A}
\newcommand{\En}{\mathcal E_d}
\newcommand{\Rel}{\operatorname{Rel}}
\title[$E_d$-Enhancements and Secondary Polytopes]{$E_d$-Enhancements of $L_\infty$-Algebras Associated to Secondary Polytopes}
\author{Longfei Li}
\date{}
\begin{document}
\begin{abstract}
We study $E_d$-enhancements of the geometric $L_\infty$-algebras
associated to secondary polytopes of finite configurations in $\R^d$,
$d\geq2$.
We construct an action of the semialgebraic chains of the Fulton--MacPherson
operad on the shifted algebra $B_A=\gA[1-d]$ when $|A|\leq3d$.
The action recovers the prescribed higher Lie operations on fundamental
chains. The construction is obtained from a general extension theorem
for weighted shifted $L_\infty$-algebras, using reflection to control
relative cohomology obstructions. For $|A|\leq3d-1$, the resulting
reflection-compatible integration actions are homotopic relative to
their shifted Lie structure. We also give an abstract enhancement
without a cardinality bound, direct formulas for circuits, and explicit
operations for convex pentagons and hexagons.
\end{abstract}
\maketitle
\tableofcontents
\section{Introduction}
\label{sec:introduction}

\subsection{Background}
The algebra of the infrared was introduced by Gaiotto--Moore--Witten
\cite{GMW} to describe algebraic structures arising from massive
two-dimensional $\mathcal N=(2,2)$ field theories. The vacua determine
a finite configuration in the complex plane, and planar webs organize
the interactions between soliton states. The resulting operations
satisfy higher Lie identities. With suitable ordering and boundary
data, the same formalism also produces directed $A_\infty$-algebras
and their deformations.

Kapranov--Kontsevich--Soibelman gave a construction of these structures
using secondary polytopes \cite{KKS}. The faces of a secondary polytope
are indexed by regular marked subdivisions. Under the factorization
hypothesis, a face is a product of the secondary polytopes of the
cells of its subdivision. Cellular incidence maps then define
multilinear operations, and the cellular boundary identity gives the
$L_\infty$ relations. Since secondary polytopes are defined for point
configurations in $\mathbb R^d$, this construction provides a
higher-dimensional version of its $L_\infty$-algebra.

The operadic background comes from the theory of iterated loop
spaces. May's little cubes operads encode families of operations
parametrized by configurations \cite{May}. The homology of little
$d$-disks, for $d\geq2$, is the $d$-Poisson operad: it has a
commutative product and a Lie bracket of cohomological degree
$1-d$ \cite[Part~III]{CLM}. Formality relates this homology operad
to chains on little disks \cite{Kontsevich,LV}. Our enhancement
problem retains the higher Lie operations themselves and asks
for compatible operations on all configuration chains.

Kapranov--Soibelman announced an $E_d$ refinement using factorization
algebras and operadic localization \cite{KS2016}. We study this
enhancement for the shifted scalar geometric $L_\infty$-algebra
associated to secondary polytopes. Our construction assigns operations
to chains in compactified configuration spaces and recovers the
prescribed higher Lie operations on their fundamental chains.
The scope of this paper is the $L_\infty$ part of the algebra of
the infrared. Enhancements of its $A_\infty$-algebras and mixed
operations are left for future work.

\subsection{The enhancement problem}
Let $d\geq2$, and let $A\subset\mathbb R^d$ be a finite,
affinely spanning configuration in affine general position.
We work over $\mathbb R$ and use the scalar geometric $L_\infty$-algebra
$\mathfrak g_A$. Its summands are indexed by full-dimensional
polytopes $Q$ with vertices in $A$, carrying the full marking
$A\cap Q$. Put
\[
 B_A=\mathfrak g_A[1-d].
\]
The differential on $B_A$ is zero. The secondary-polytope operations
give maps
\[
 \tau_r:B_A^{\otimes r}\longrightarrow B_A,
 \qquad |\tau_r|=1-d(r-1),\qquad r\geq2,
\]
which define an action of the shifted homotopy Lie operad
$\mathsf L_{\infty,d}$.

We impose the following factorization hypothesis throughout the
secondary-polytope construction. For every full-dimensional
$A'\subseteq A$ and every regular marked subdivision
$\{(Q_i,A_i)\}$ of $(\operatorname{Conv}(A'),A')$, the height
restriction map
\[
 \mathbb R^{A'}/\operatorname{Aff}(A')
 \longrightarrow
 \bigoplus_i\mathbb R^{A_i}/\operatorname{Aff}(A_i)
\]
is surjective. Here $\operatorname{Aff}$ denotes restrictions of
affine-linear functions. This condition gives compatible product
identifications of secondary-polytope faces. The correction in
\cite[Proposition~5.2.14]{KSS} is stated for planar configurations.
In Section~\ref{sec:factorization-criterion} we give the criterion
in arbitrary dimension: the kernel of the sum of the cell relation
spaces measures exactly the failure of the face dimensions to add.
Its vanishing is equivalent to the canonical affine product
decomposition. In dimension two this defect agrees with the
exceptionality of \cite[Section~5.2]{KSS}.
It suffices to test subdivisions of $A$ itself: the condition passes
to every full-dimensional subconfiguration by
Corollary~\ref{cor:ambient-factorization-test}. Thus for $d=2$
the hypothesis is exactly exceptionally general position in the
sense of \cite[Definition~5.2.11]{KSS}.

For $r\geq2$, let
\[
 M_r=\operatorname{FM}_d(r),\qquad D_r=d(r-1)-1,
\]
where $\operatorname{FM}_d$ is the Fulton--MacPherson model of the
little $d$-disks operad. Its boundary strata record configurations
with nested clusters. We use the dg operad
\[
 \mathcal E_d(r)=C^{\mathrm{SA}}_{-*}(M_r;\mathbb R)
\]
of semialgebraic geometric chains, with a $p$-chain in cohomological
degree $-p$. The unary space is a point, and all operads in this
paper have positive arities. The standard fundamental-chain morphism
\cite[Section~2.2]{CW} is
\[
 \iota_d:\mathsf L_{\infty,d}\longrightarrow\mathcal E_d,
 \qquad \mu_r\longmapsto[M_r].
\]
Our question is how to construct an action of $\mathcal E_d$ on
$B_A$ whose restriction along $\iota_d$ is the given action $\tau$.

\subsection{Main algebraic results}
The first result constructs the operations directly on configuration
chains.

\begin{theorem}[Geometric enhancement]
\label{thm:intro-geometric}
Let $A\subset\mathbb R^d$, $d\geq2$, be finite, affinely spanning,
and in affine general position. Assume the factorization hypothesis
above. If $|A|\leq3d$, there is a dg-operad morphism
\[
 \rho_A:\mathcal E_d\longrightarrow\operatorname{End}_{B_A}
\]
such that
\[
 \rho_{A,r}([M_r])=\tau_r\qquad(r\geq2).
\]
The action is obtained by integrating smooth endomorphism-valued
forms. Its degree-zero binary product is zero, and the unary
identity acts as $\operatorname{id}_{B_A}$.
\end{theorem}

For $|A|\leq d+3$, the factorization hypothesis follows from affine
general position. Thus every affinely spanning configuration in
this smaller range admits the enhancement of
Theorem~\ref{thm:intro-geometric}. For vertices of a convex polygon, factorization follows by
successively assembling heights across the diagonals of a dissection.
In dimension two the theorem therefore includes convex hexagons.

There is also an algebraic construction without a cardinality bound.
Let $P_d^\infty$ be the nonunital homotopy $d$-Poisson operad,
with its standard inclusion
$i:\mathsf L_{\infty,d}\to P_d^\infty$. The cobar models of
\cite[Section~2]{TW} give a retraction
$p:P_d^\infty\to\mathsf L_{\infty,d}$ with $pi=\operatorname{id}$.
Real operadic formality identifies $P_d^\infty$ as a dg model for
$E_d$ \cite{LV}.

\begin{proposition}[Abstract enhancement]
\label{prop:intro-algebraic}
Let $A$ satisfy the affine and factorization assumptions above,
with arbitrary cardinality. The composite
\[
 \rho_A^{\mathrm{alg}}=\tau p:
 P_d^\infty\longrightarrow\operatorname{End}_{B_A}
\]
defines a nonunital $E_d$-enhancement in this operad model. Its
restriction along $i$ is exactly $\tau$, and its binary product
is zero.
\end{proposition}

The proposition follows from the standard operad models. The
geometric theorem gives an action on the specified configuration
chains, with the values of the fundamental chains fixed in advance.

The relation between Lie and $E_d$-algebras also has a universal
construction. Knudsen constructs higher enveloping algebras as a
left adjoint to the shifted Lie functor, with a description in
terms of Lie algebra homology \cite{Knudsen}. Here the enhancement
is required on the specified complex $B_A$, with its primary
operations fixed. A comparison with the higher enveloping
construction is a separate question, left for future work.

\subsection{The geometric construction and examples}
The construction uses the dimension of a secondary polytope as a
weight. For a geometric marked polytope $Q$, set
\[
 w(Q)=|A\cap Q|-d-1.
\]
The associated generator of $B_A$ has degree $d+w(Q)$, and every
primary operation satisfies
\[
 w_{\mathrm{out}}=1+\sum_i w_{\mathrm{in},i}.
\]
The weights therefore lie between zero and $W=|A|-d-1$.

More generally, the argument applies to a finite-dimensional graded
real vector space $V=\bigoplus_{a=0}^W V_{(a)}$ with zero
differential and a fixed $\mathsf L_{\infty,d}$-action whose primary
operations increase total weight by one. A weight-$w$ operation
in arity $r$ is represented by a closed form $\alpha_{r,w}$ of
degree $d(r-1)-w$, with values in endomorphisms of the opposite
cohomological degree and weight increase $w$. At weight one we take
\[
 \alpha_{r,1}=\eta_r\tau_r,\qquad \int_{M_r}\eta_r=1,
\]
where $\eta_r$ is an equivariant top form supported in the interior.
At higher weights, compositions of the previously constructed forms
prescribe the boundary values of $\alpha_{r,w}$. Closed extensions
produce the remaining chain operations by integration.

At weight two, the boundary integral is the quadratic
$L_\infty$ identity. At subsequent weights, we use a fixed
hyperplane reflection $s$ of the auxiliary $\mathbb R^d$ and require
\[
 s^*\alpha_{r,w}=(-1)^{r-1}\alpha_{r,w},
\]
with trivial reflection action on the endomorphism coefficients.
Poincare--Lefschetz duality places the extension obstruction in the
reflection-invariant part of $H_{w-2}(M_r;\mathbb R)$, with the
corresponding endomorphism coefficients. This invariant homology
vanishes for $3\leq w\leq2d-1$. The construction therefore gives
an $\mathcal E_d$-action whenever $W\leq2d-1$, yielding
Theorem~\ref{thm:intro-geometric}.

The same method compares the choices of forms by extending over
$M_r\times[0,1]$.

\begin{proposition}[Comparison of geometric enhancements]
\label{prop:intro-comparison}
Fix $A$ satisfying the affine and factorization assumptions, with
$|A|\leq3d-1$, and fix its shifted Lie action $\tau$ and the
auxiliary reflection $s$. Any two integration actions constructed
from normalized primary forms and closed correction forms of the
prescribed weights, satisfying the permutation, insertion, and
reflection conditions above, are polynomially homotopic relative
to $\tau$. More precisely, there is a dg-operad morphism
\[
 \mathcal H:\mathcal E_d\longrightarrow
 \operatorname{End}_{B_A}\otimes\mathbb R[t,dt]
\]
whose endpoint evaluations are the two actions and for which
$\mathcal H\iota_d=\tau\otimes1$. Here $|t|=0$, $|dt|=1$,
and $d(t)=dt$.
\end{proposition}

The comparison is proved in Theorem~\ref{thm:reflected-comparison}.
It allows the normalized primary forms to vary.
The relative cylinder calculation replaces $H_{w-2}(M_r)$ by
$H_{w-1}(M_r)$, which gives the bound $W\leq2d-2$.

For an affine circuit, the secondary polytope is an interval.
Every primary operation lands in the line corresponding to the
whole circuit, and every operation with that line as an input is
zero. Consequently, the action has the direct formula
\[
 \rho_r(c)=
 \begin{cases}
   \bigl(\int_c\eta_r\bigr)\tau_r,&\dim c=D_r,\\
   0,&\dim c\neq D_r,
 \end{cases}
 \qquad r\geq2.
\]
Together with the unary identity, this defines the entire action.

A triangle with one interior point gives an $E_2$-algebra with
zero binary product and bracket, but a nonzero ternary
fundamental-chain operation. Its ternary Lie operation persists
under quasi-isomorphisms between minimal $L_\infty$-algebras.
A triangular bipyramid gives an $E_3$-example with both binary and
ternary primary operations.

The convex pentagon exhibits the first boundary correction. Two
successive binary gluings have a nonzero composite, so the
corresponding two-dimensional boundary chain in
$\operatorname{FM}_2(3)$ must act nontrivially. This operation is
realized by a closed two-form, while the ternary fundamental-chain
operation is zero. The convex hexagon then illustrates the next
weight of the reflected extension construction.

\subsection{Factorization algebras and possible applications}
\label{sec:intro-related-contexts}
Factorization algebras provide a local formulation of the structures
considered here. They assign complexes to open sets, together with
compatible multiplication maps for disjoint open subsets, and satisfy
a descent condition. Local constancy requires the maps associated to
inclusions of disks to be quasi-isomorphisms. Costello--Gwilliam
develop this language and its relation to observables in field theory
\cite[Chapter~6]{CG}. The comparison with little disks is expressed
by Lurie's equivalence between disk algebras and locally constant
disk data \cite[Theorem~5.4.5.9]{HA}; its nonunital version is
\cite[Proposition~5.4.5.15]{HA}. Thus the enhancements studied here
also have a local interpretation, with the same positive-arity
convention as our operads.

Factorization homology supplies a corresponding construction on
manifolds. For a unital $E_d$-algebra $C$ and a framed $d$-manifold
$M$, it forms $\int_M C$ by assembling the operations on disks;
its excision property describes the effect of gluing manifolds
\cite[Section~3]{AF}. In dimension one, its value on the circle
is the Hochschild homology of the coefficient algebra. These results
suggest using the secondary-polytope algebras as coefficients in
manifold invariants. For our nonunital algebras, such an application
must retain the nonunital convention or pass to an augmented unital
model. The computation of these invariants is left for future work.
We use framed manifolds in this discussion: the chosen reflection
in the extension argument is not an action of the full orthogonal
group on the algebra.

There is a second connection through deformation theory. The
Hochschild cochains of an associative algebra admit an $E_2$-structure,
as in the proof of the Deligne conjecture by Kontsevich--Soibelman
\cite{KSDeligne}. More generally, higher Hochschild cochains of a
unital $E_d$-algebra carry an $E_{d+1}$-structure, and Francis relates
them to the tangent complex of the algebra \cite{Francis}.
The explicit operations constructed in this paper give concrete
examples for studying these deformation-theoretic constructions,
using an augmented unital model. Their calculation, and
comparison with deformation problems arising from the original
point configuration, remain future work.

Finally, Voronov's Swiss-cheese operad combines the little disks
and little intervals operads, describing operations with inputs
in the interior of the upper half-plane and on its boundary
\cite{Voronov}. This is relevant to the $A_\infty$ and mixed parts
of the algebra of the infrared \cite{KKS,KSS}. Constructing a
compatible action of such a two-colored operad requires additional
boundary and mixed operations. We leave that extension to future
work.

\subsection{Outline of this paper}
Section~\ref{sec:secondary-polytopes} recalls secondary polytopes, the factorization hypothesis,
and the full and geometric $L_\infty$-algebras. It also fixes the
grading and orientation conventions.

Section~\ref{sec:operad-models} introduces the operad models, the fundamental-chain
morphism, and the abstract enhancement. Section~\ref{sec:configuration-space-action} gives the direct
integration construction for circuits.

Section~\ref{sec:weight-extension} develops the weight filtration, the relative extension
problem, and the reflection argument proving the geometric
enhancement theorem. Section~\ref{sec:geometric-homotopies} studies homotopies between the
constructed actions and proves the comparison proposition.

Section~\ref{sec:examples} presents the explicit configurations and their operations.
Section~\ref{sec:further-directions} describes the subsequent extension problem and the
directions for further work.

\subsection*{Conventions}
All vector spaces and chains have real coefficients. We use
cohomological grading, with $(V[s])^j=V^{j+s}$. Tensor products,
permutations, and compositions of homogeneous maps include Koszul
signs. Operads retain their unary identity and have no nullary
operations. The dimension $d$ is at least two.

\subsection*{Acknowledgements}
I am grateful to my advisor, Yan Soibelman, for introducing me to this topic and for his helpful comments. I used GPT to assist with revising parts of the manuscript and correcting typographical errors.
\section{Secondary polytopes and the geometric
\texorpdfstring{$L_\infty$}{L-infinity}-algebra}
\label{sec:secondary-polytopes}

We recall the combinatorial construction needed for the
enhancement theorem. We use the secondary-polytope
description of \cite{GKZ,DLRS} and the algebraic framework of
\cite{KKS}. We state a factorization hypothesis explicitly,
prove that it holds for the small configurations considered
here, and fix the duality and suspension conventions.

\subsection{Marked subdivisions}
Fix $d\geq2$ and a finite, affinely spanning configuration
$A\subset\R^d$. Unless a statement explicitly relaxes this assumption, $A$ is in
\emph{affine general position}: every subset of at most
$d+1$ points is affinely independent.

\begin{definition}\label{def:marked-polytope}
A \emph{marked subpolytope} is a pair $(Q_B,B)$, where
$B\subseteq A$, $Q_B=\Conv(B)$, and $\dim Q_B=d$.
It is \emph{geometric relative to $A$} if
$B=A\cap Q_B$. A full-dimensional simplex $\Delta$
with vertices in $A$ is \emph{empty} if
$A\cap\Delta=\operatorname{Vert}(\Delta)$.
\end{definition}

A \emph{marked subdivision} of $(Q_B,B)$ is a collection
$\mathcal S=\{(Q_i,B_i)\}_{i\in I}$ of full-dimensional
marked subpolytopes, with $B_i\subseteq B$, such that
the $Q_i$ cover $Q_B$, distinct cells meet in a common
proper face, possibly empty, and
\[
B_i\cap(Q_i\cap Q_j)=B_j\cap(Q_i\cap Q_j).
\]
We do not require $\bigcup_iB_i=B$. A subdivision with
one maximal cell and a smaller marking is allowed.
A \emph{triangulation} is a marked subdivision in
which every $B_i$ consists of $d+1$ points.

For a height function $h:B\to\R$, define
\[
\underline h(x)=
\min_{\substack{t_b\geq0,\ \sum_b t_b=1\\
                 \sum_b t_b b=x}}
\sum_{b\in B}t_bh(b).
\]
Its maximal domains of affinity are convex polytopes
$Q_i$. Mark each of them by the contact set
\[
B_i=\{b\in B\cap Q_i:h(b)=\underline h(b)\}.
\]
A subdivision obtained in this way is \emph{regular}.
The markings record which points lie on the lower
faces of the lifted configuration; a point lying
strictly above the lower envelope is omitted. For the lifting
construction and the secondary fan, see \cite[Chapter~5]{DLRS}.

We write $\mathcal T\preceq\mathcal S$ when
$\mathcal T$ refines $\mathcal S$, including its
markings. Thus the cells of $\mathcal T$ inside
$(Q_i,B_i)$ form a marked subdivision using subsets
of $B_i$. The trivial subdivision $\{(Q_B,B)\}$
is the largest element. A regular subdivision is
\emph{coarse} if it is maximal among the proper
regular subdivisions.

For a triangulation $\mathcal T$, its GKZ vector is
\[
\phi_{\mathcal T}(b)=
\sum_{\substack{\Delta\in\mathcal T\\
                 b\in\operatorname{Vert}(\Delta)}}
\operatorname{Vol}_d(\Delta).
\]
The secondary polytope is
\[
\Sigma(B)=
\Conv\{\phi_{\mathcal T}:\mathcal T
       \text{ a regular triangulation of }(Q_B,B)\}
\subset\R^B.
\]
An unused mark has GKZ coordinate zero.

We use the following standard properties of secondary
polytopes \cite[Chapter~7]{GKZ} and \cite[Chapter~5]{DLRS}:
\begin{equation}\label{eq:secondary-dimension}
w(B):=\dim\Sigma(B)=|B|-d-1,
\end{equation}
and regular marked subdivisions correspond to
nonempty faces, with
\[
F_{\mathcal T}\subseteq F_{\mathcal S}
\quad\Longleftrightarrow\quad
\mathcal T\preceq\mathcal S.
\]
In particular, triangulations give vertices and
coarse subdivisions give facets. These statements
use the specified refinement order. The number
of maximal cells alone does not determine the
dimension of the corresponding face.

\subsection{Factorization of faces}
Write $\operatorname{Aff}(B)\subset\R^B$ for the
restrictions of affine-linear functions on $\R^d$.

\begin{definition}[Factorization hypothesis]
\label{def:factorization}
For every full-dimensional $B\subseteq A$ and every
regular marked subdivision
$\mathcal S=\{(Q_i,B_i)\}_{i\in I}$ of $(Q_B,B)$,
we require the restriction map
\begin{equation}\label{eq:height-restriction}
R_{\mathcal S}:
\R^B/\operatorname{Aff}(B)
\longrightarrow
\bigoplus_{i\in I}\R^{B_i}/\operatorname{Aff}(B_i)
\end{equation}
to be surjective.
\end{definition}

\begin{proposition}\label{prop:face-factorization}
Under Definition~\ref{def:factorization}, extension
by zero followed by addition gives an affine
isomorphism
\[
\prod_{i\in I}\Sigma(B_i)\xrightarrow{\ \sim\ }
F_{\mathcal S}.
\]
These isomorphisms are compatible with further
regular refinements. Consequently,
\begin{equation}\label{eq:face-dimension}
\dim F_{\mathcal S}=\sum_{i\in I}w(B_i).
\end{equation}
\end{proposition}

\begin{proof}
Proposition~\ref{prop:relation-criterion} below proves the affine
isomorphism for each regular subdivision whose height map is
surjective. Definition~\ref{def:factorization} supplies this
condition for every marking and subdivision under consideration.
Under further refinements, extension by zero and addition are
associative, giving the compatibility statement.
\end{proof}

\begin{remark}\label{rem:factorization-hypothesis}
Affine general position alone does not imply this hypothesis.
Section~\ref{sec:factorization-criterion} gives the exact
dimension-independent condition and its relation to the planar
correction of \cite[Proposition~5.2.14]{KSS}. The small-cardinality
case follows from the elementary count below.
\end{remark}

\begin{lemma}\label{lem:small-factorization}
If $|A|\leq d+3$, then $A$ satisfies
Definition~\ref{def:factorization}.
\end{lemma}

\begin{proof}
Consider a marked subdivision of a full-dimensional
$B\subseteq A$. Two distinct maximal cells meet
in a proper face, whose affine span has dimension
at most $d-1$. General position implies
$|B_i\cap B_j|\leq d$. If both $|B_i|\geq d+2$
and $|B_j|\geq d+2$, then
\[
|B_i\cup B_j|\geq 2(d+2)-d=d+4,
\]
contrary to $|B|\leq d+3$.
Thus at most one cell has more than $d+1$ marks.
Every other cell is a simplex and has
$\R^{B_i}/\operatorname{Aff}(B_i)=0$.
If a nonzero summand remains in the target of
\eqref{eq:height-restriction}, restriction onto
that summand is surjective: extend any height
function on its marking to $B$. This also
includes the trivial subdivision.
\end{proof}

\subsection{The exact factorization criterion in arbitrary dimension}
\label{sec:factorization-criterion}
We now identify the precise defect of the product decomposition.
This gives a criterion for an individual regular subdivision, in
every ambient dimension. It also relates the height condition to
the planar correction in \cite[Proposition~5.2.14]{KSS}.
For a full-dimensional marking $B$, let
\[
 \Rel(B)=\left\{(c_b)_{b\in B}:\sum_{b\in B}c_b=0,
                       \quad\sum_{b\in B}c_b b=0\right\}.
\]
Extension by zero identifies $\Rel(B_i)$ with a subspace of
$\Rel(B)$ whenever $B_i\subseteq B$.
For a regular marked subdivision $\mathcal S=\{(Q_i,B_i)\}$, set
\begin{equation}\label{eq:factorization-defect}
\begin{aligned}
 J_{\mathcal S}:\bigoplus_i\Rel(B_i)&\longrightarrow\Rel(B),
       &(c_i)_i&\longmapsto\sum_i c_i,\\
 e(\mathcal S)&=\dim\ker J_{\mathcal S}.
\end{aligned}
\end{equation}
We call $e(\mathcal S)$ the \emph{factorization defect}.

\begin{proposition}[Factorization criterion]\label{prop:relation-criterion}
Let $B\subset\R^d$ be finite and affinely spanning, and let
$\mathcal S$ be a regular marked subdivision of $(Q_B,B)$.
Then
\begin{equation}\label{eq:factorization-dimension-defect}
 \dim F_{\mathcal S}
 =\sum_i\bigl(|B_i|-d-1\bigr)-e(\mathcal S).
\end{equation}
The following conditions are equivalent:
\begin{enumerate}
\item $R_{\mathcal S}$ is surjective;
\item $J_{\mathcal S}$ is injective, or equivalently $e(\mathcal S)=0$;
\item $\dim F_{\mathcal S}=\sum_i\dim\Sigma(B_i)$;
\item extension by zero followed by addition is an affine isomorphism
      $\prod_i\Sigma(B_i)\longrightarrow F_{\mathcal S}$;
\item $F_{\mathcal S}$ is combinatorially equivalent to
      $\prod_i\Sigma(B_i)$.
\end{enumerate}
These assertions do not require affine general position.
\end{proposition}

\begin{proof}
The annihilator of $\operatorname{Aff}(B)$ in $(\R^B)^*$ is
$\Rel(B)$. Under this identification, the dual of restriction
of heights is extension by zero of relations. Thus
$J_{\mathcal S}=R_{\mathcal S}^*$, proving the equivalence of
(1) and (2).

To calculate the dimension of the face, put
\[
 D_{\mathcal S}
 =\{h\in\R^B:h|_{B_i}\in\operatorname{Aff}(B_i)\text{ for every }i\}.
\]
The secondary cone consisting of heights inducing exactly
$\mathcal S$ is a nonempty relatively open subset of $D_{\mathcal S}$.
Indeed, for each $i$ the affine function interpolating $h|_{B_i}$
is uniquely determined. Heights inducing $\mathcal S$ are characterized,
within $D_{\mathcal S}$, by the finitely many strict inequalities
placing every other lifted mark above each such supporting affine
function. Regularity gives a solution. Matching heights at the
vertices of common faces ensures that the interpolating functions
agree on those faces. Marks omitted from all the $B_i$ have no
linear equality imposed on their heights.

The secondary fan is the normal fan of $\Sigma(B)$
\cite[Chapter~7]{GKZ}. Consequently, the linear span of this
normal cone is $D_{\mathcal S}$, and the direction space of its
dual face is
\[
 \operatorname{dir}(F_{\mathcal S})
 =D_{\mathcal S}^{\perp}
 =\sum_i\Rel(B_i)=\operatorname{im}J_{\mathcal S}.
\]
Here the annihilator and the relation spaces are identified inside
$\R^B$ by the coordinate pairing. Since every cell is full-dimensional,
$\dim\Rel(B_i)=|B_i|-d-1=\dim\Sigma(B_i)$.
Rank--nullity proves \eqref{eq:factorization-dimension-defect}
and the equivalence of (2) and (3).

For clarity, (1) gives the asserted affine isomorphism directly.
Prescribe a regular triangulation of each cell and choose heights
realizing it. Surjectivity supplies a height function on $B$ whose
restrictions differ from these chosen heights by affine functions.
A sufficiently small perturbation of a height function for
$\mathcal S$ then realizes all the chosen triangulations
simultaneously. This is the usual refinement property of lower
faces: the strict inequalities separating distinct original cells
persist, while the perturbation subdivides each original lower face.
Thus every tuple of vertices of $\prod_i\Sigma(B_i)$ maps to the
GKZ vector of a regular triangulation refining $\mathcal S$.
Conversely, each vertex of $F_{\mathcal S}$ restricts to such a
tuple, and its GKZ vector is the sum of the cell GKZ vectors.
The addition map therefore has image $F_{\mathcal S}$.
Its linear part on direction spaces is $J_{\mathcal S}$, which
is injective. This proves (4). An affine isomorphism implies (5),
and combinatorial equivalence preserves dimension, so (5)
implies (3). This proves all the equivalences.
\end{proof}

The affine identifications are compatible with iterated refinements
whenever the corresponding defects vanish: both ways of assembling
the GKZ vectors extend by zero and add the same cell contributions.
Thus Definition~\ref{def:factorization} is exactly the vanishing of
$e(\mathcal S)$ for the regular subdivisions of every
full-dimensional subconfiguration used in the algebra.

\subsubsection{Restriction to subconfigurations}
The quantification over subconfigurations in
Definition~\ref{def:factorization} can be reduced to subdivisions
of the ambient configuration.

\begin{lemma}[Extending a regular subdivision]
\label{lem:extend-subconfiguration}
Let $B\subseteq A\subset\R^d$ be finite, with $B$ affinely
spanning. Every regular marked subdivision $\mathcal S$ of
$(\Conv(B),B)$ occurs as a collection of maximal cells of a regular
marked subdivision $\mathcal T$ of $(\Conv(A),A)$. Moreover,
\[
 \ker J_{\mathcal S}\lhook\joinrel\longrightarrow\ker J_{\mathcal T},
 \qquad e(\mathcal S)\leq e(\mathcal T).
\]
\end{lemma}

\begin{proof}
Choose heights $h:B\to\R$ inducing $\mathcal S$. For its maximal
cells $(Q_i,B_i)$, let $\lambda_i:\R^d\to\R$ be the affine
supporting functions of the lower hull. Thus
$h(b)\geq\lambda_i(b)$ for every $b\in B$, with equality exactly
on $B_i$. Extend $h$ to $A$ by choosing
\[
 \widehat h(a)>\max_i\lambda_i(a)\qquad(a\in A\setminus B).
\]
Each $\lambda_i$ remains a lower supporting function with exactly
the same contact set $B_i$. Since $Q_i$ is full-dimensional,
$(Q_i,B_i)$ is a maximal cell of the regular subdivision
$\mathcal T$ defined by $\widehat h$.

A tuple of cell relations in $\ker J_{\mathcal S}$ extends by
zero on the additional cells of $\mathcal T$. The same relation
vectors, extended by zero from $B$ to $A$, still sum to zero.
This gives the asserted injective map of kernels and the inequality.
\end{proof}

\begin{corollary}[A test on the ambient configuration]
\label{cor:ambient-factorization-test}
A finite affinely spanning configuration $A$ satisfies
Definition~\ref{def:factorization} if and only if
$e(\mathcal T)=0$ for every regular marked subdivision
$\mathcal T$ of $(\Conv(A),A)$.
\end{corollary}

\begin{proof}
The forward implication is Proposition~\ref{prop:relation-criterion}
with $B=A$. Conversely, for any full-dimensional $B\subseteq A$
and any regular subdivision $\mathcal S$ of $B$, choose
$\mathcal T$ as in Lemma~\ref{lem:extend-subconfiguration}.
Then $0\leq e(\mathcal S)\leq e(\mathcal T)=0$.
Proposition~\ref{prop:relation-criterion} gives the required
surjectivity for every $B$ and $\mathcal S$.
\end{proof}

The criterion is a finite rank computation once the regular
subdivisions are known. For a cell marking $B_i$, form the
$(d+1)\times|B_i|$ matrix with columns $(1,b)^{\top}$,
$b\in B_i$. Choose a basis of its kernel, extend these vectors
by zero to $\R^A$, and concatenate the resulting columns over
the cells. If $N_{\mathcal T}$ is this matrix, then
\[
 e(\mathcal T)
 =\sum_i(|B_i|-d-1)-\operatorname{rank}N_{\mathcal T}.
\]
This also gives an exact arithmetic test for rational configurations.

\subsubsection{Comparison with the planar correction}
The corrected statement in \cite[Proposition~5.2.14]{KSS} is
formulated for $A\subset\mathbb C=\R^2$ in exceptionally general
position; its conclusion is a combinatorial product decomposition.
The following calculation identifies its subdivision-wise condition
with the defect above. It also accounts for unused marks explicitly.

\begin{proposition}\label{prop:planar-exceptionality}
For a planar regular marked subdivision $\mathcal S$ in affine
general position,
\[
 e(\mathcal S)=\operatorname{exc}_{\mathrm{KSS}}(\mathcal S),
\]
where $\operatorname{exc}_{\mathrm{KSS}}(\mathcal S)$ is the
dimension of the first cohomology of the deformation complex
of \cite[Section~5.2.B]{KSS}.
\end{proposition}

\begin{proof}
Write $U=\bigcup_iB_i$. Let $f$ be the number of maximal polygons,
$a$ the number of internal edges, and $v$ the number of interior
vertices of the underlying subdivision. Let
$\operatorname{PL}(\mathcal S)$ be the space of continuous functions
on $Q_B$ that are affine on each cell. Restriction of a height in
$D_{\mathcal S}$ to $U$ determines such a function, whereas its
values on $B\setminus U$ are free. Hence
\[
 \dim D_{\mathcal S}
 =\dim\operatorname{PL}(\mathcal S)+|B\setminus U|,
 \qquad
 \dim F_{\mathcal S}=|U|-\dim\operatorname{PL}(\mathcal S).
\]
The Euler characteristic calculation for the planar deformation
complex, whose terms have dimensions $3f,2a,v$ and whose second
cohomology vanishes, gives
\[
 \dim\operatorname{PL}(\mathcal S)
 =3f-2a+v+\operatorname{exc}_{\mathrm{KSS}}(\mathcal S)
\]
\cite[Propositions~5.2.2 and~5.2.4]{KSS}, applied after discarding
unused marks. Affine general position implies that a mark on an
edge is one of its endpoints. Counting cell vertices along edges,
and counting each marked point in a cell interior once, therefore gives
\[
 \sum_i|B_i|=|U|+2a-v.
\]
Combining these equations with
\eqref{eq:factorization-dimension-defect} proves the assertion.
\end{proof}

In particular, the planar non-exceptionality condition yields the
canonical affine product map in Proposition~\ref{prop:relation-criterion}.
For arbitrary $d$, the relation map $J_{\mathcal S}$ supplies the
corresponding criterion directly.

\begin{corollary}\label{cor:kss-factorization-equivalence}
For a finite affinely spanning planar configuration $A$ in affine
general position, Definition~\ref{def:factorization} is equivalent
to exceptionally general position in the sense of
\cite[Definition~5.2.11]{KSS}.
\end{corollary}

\begin{proof}
The latter condition is the vanishing of
$\operatorname{exc}_{\mathrm{KSS}}(\mathcal T)$ for every regular
subdivision of $A$. Apply Proposition~\ref{prop:planar-exceptionality}
and Corollary~\ref{cor:ambient-factorization-test}.
\end{proof}

\subsubsection{A failure of the criterion}
\begin{example}[A six-point failure of factorization]
\label{ex:factorization-failure}
Consider
\[
\begin{aligned}
 a&=(3,0),& b&=(0,3),& c&=(-3,-3),\\
 a'&=(1,0),& b'&=(0,1),& c'&=(-1,-1).
\end{aligned}
\]
These points are in affine general position. One direct check uses
the three vectors $u=(1,0)$, $v=(0,1)$, $z=(-1,-1)$: the points
are $u,v,z,3u,3v,3z$. A triple with two points on one of these rays
has nonzero determinant since its third point is on a different
ray. A triple with one point on each ray has the form
$\alpha u,\beta v,\gamma z$, with positive $\alpha,\beta,\gamma$;
its affine determinant is
$\alpha\beta+\alpha\gamma+\beta\gamma>0$.

Assign height one to $a,b,c$ and zero to $a',b',c'$.
The lower envelope on $\Conv\{a,b,c\}$ is
\[
 \max\left\{0,\frac{x+y-1}{2},
               \frac{-2x+y-1}{2},\frac{x-2y-1}{2}\right\}.
\]
Indeed, each displayed affine function lies below the assigned
heights, and its contact points span the corresponding cell.
The resulting regular subdivision consists of the inner triangle
and the three quadrilaterals
\[
 abb'a',\qquad bcc'b',\qquad caa'c'.
\]
In the coordinate order $(a,b,c,a',b',c')$, their relation spaces
are generated by
\[
\begin{aligned}
 c_{ab}&=(1,-1,0,-3,3,0),\\
 c_{bc}&=(0,1,-1,0,-3,3),\\
 c_{ca}&=(-1,0,1,3,0,-3).
\end{aligned}
\]
Their sum is zero. The direct sum of the three one-dimensional
relation spaces therefore maps noninjectively to $\Rel(A)$, so
$R_{\mathcal S}$ is not surjective. The first two displayed relations
are independent, so $e(\mathcal S)=1$ and
$\dim F_{\mathcal S}=3-1=2$. The product of the three secondary
intervals has dimension three, and thus cannot be combinatorially
equivalent to this face. This is the concentric-triangle
example of \cite[Example~5.2.10]{KSS}, written with rational coordinates.
\end{example}

\begin{figure}[htbp]
\centering
\begin{tikzpicture}[scale=.62,every node/.style={font=\small}]
 \coordinate (a) at (3,0); \coordinate (b) at (0,3);
 \coordinate (c) at (-3,-3); \coordinate (ap) at (1,0);
 \coordinate (bp) at (0,1); \coordinate (cp) at (-1,-1);
 \fill[gray!15] (ap)--(bp)--(cp)--cycle;
 \draw[thick] (a)--(b)--(c)--cycle;
 \draw (ap)--(bp)--(cp)--cycle;
 \draw (a)--(ap) (b)--(bp) (c)--(cp);
 \foreach \v in {a,b,c,ap,bp,cp} \fill (\v) circle (2pt);
 \node[right] at (a) {$a$}; \node[above] at (b) {$b$};
 \node[below left] at (c) {$c$}; \node[below right] at (ap) {$a'$};
 \node[above left] at (bp) {$b'$}; \node[below right] at (cp) {$c'$};
\end{tikzpicture}
\caption{A regular subdivision for which the three cell relation
spaces are linearly dependent. Ordinary affine general position
does not imply the height-surjectivity condition.}
\label{fig:factorization-failure}
\end{figure}
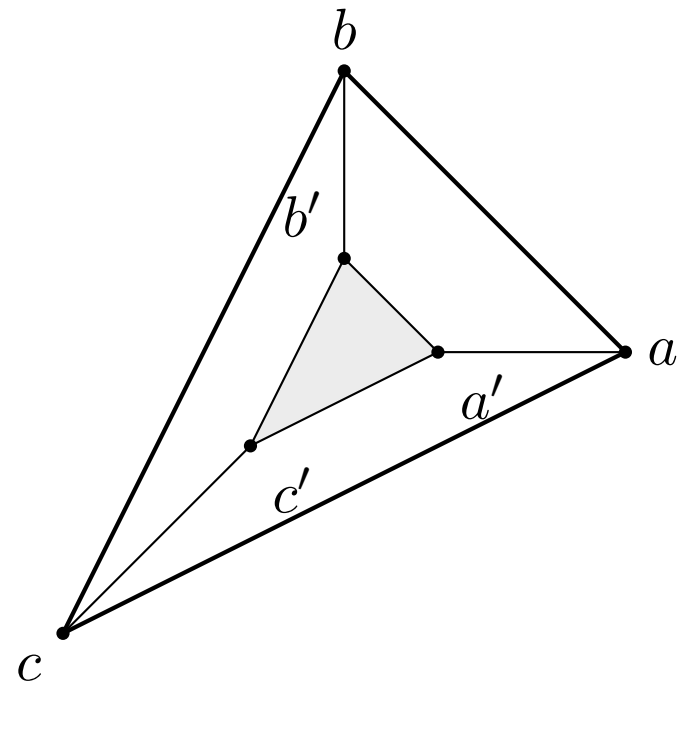

\begin{lemma}[Convex polygon configurations]\label{lem:polygon-factorization}
The vertices of a strictly convex polygon in $\R^2$ satisfy
Definition~\ref{def:factorization}.
\end{lemma}

\begin{proof}
Every full-dimensional subconfiguration is again the vertex set of
a convex polygon. A subdivision of such a polygon uses noncrossing
diagonals. Its dual graph, with vertices the maximal cells and edges
the internal diagonals, is a tree. Prescribe a height function on
each cell, and choose a root of this tree. Proceed away from the root.
At each step the next cell meets the part on the root side along its
parent diagonal; its only vertices on that side are the endpoints
of the diagonal. Subtract an affine function from the new cell's
heights so that they agree with the previously assigned heights at
those two endpoints. This is possible because the endpoints are
distinct. The resulting assignments give one height function on
the whole polygon. Its restrictions have the prescribed classes
modulo affine functions. This proves surjectivity for every
subdivision and every full-dimensional subconfiguration.
\end{proof}

\subsection{The cellular differential}
Assume Definition~\ref{def:factorization}.
For a polytope $P$ of dimension $m$, put
\[
\operatorname{or}(P)=H_m(P,\partial P;\R),
\]
initially in degree zero. The cellular boundary
induces maps between the orientation lines of
a cell and its facets. Our convention is the
outward-normal-first boundary orientation.

For each full-dimensional $B\subseteq A$, define
\[
W_B=\operatorname{or}(\Sigma(B))[w(B)],
\qquad
W_A=\bigoplus_{\substack{B\subseteq A\\\dim Q_B=d}}W_B.
\]
Thus $W_B$ is a line in degree $-w(B)$. Choose
oriented fundamental generators $v_B$.
Factorization of a facet, together with its
incidence map, defines
\begin{equation}\label{eq:cellular-derivation}
d(v_B)=
\sum_{\substack{\mathcal S=\{(Q_i,B_i)\}\\
                 \mathcal S\text{ coarse}}}
\varepsilon(\mathcal S)\prod_i v_{B_i}
\quad\text{in }\operatorname{Sym}(W_A).
\end{equation}
Here $\varepsilon(\mathcal S)$ is the incidence
sign relative to the chosen product orientation.
For a coarse subdivision,
\begin{equation}\label{eq:facet-weight}
w(B)=1+\sum_iw(B_i),
\end{equation}
so $d$ has degree $+1$.
Extend $d$ by the graded Leibniz rule and set $d(1)=0$.
Interchanging two product factors contributes
$(-1)^{w(B_i)w(B_j)}$, both for cellular
orientations and for graded multiplication.
Thus \eqref{eq:cellular-derivation} is independent
of the order chosen for its cells.

\begin{proposition}\label{prop:cellular-square}
The derivation $d$ satisfies $d^2=0$.
\end{proposition}

\begin{proof}
On a generator $v_B$, a term in $d^2(v_B)$
chooses a facet of $\Sigma(B)$ and then a
facet in one of its product factors.
Proposition~\ref{prop:face-factorization}
identifies this with a flag through a
codimension-two face of $\Sigma(B)$.
Such a face lies in exactly two facets.
The resulting incidence signs are opposite,
and the final monomial is the same.
The product boundary signs are precisely
the signs of the graded Leibniz rule.
Every coefficient therefore cancels.
The square of a degree-one derivation is a
derivation, so the result on generators
proves $d^2=0$.
\end{proof}

The sum in \eqref{eq:cellular-derivation} includes
one-cell subdivisions. For example, if
$b\in B$ lies in the interior of $Q_B$, the
one-cell subdivision $(Q_B,B\setminus\{b\})$
is obtained by assigning height one to $b$
and height zero to the other marks.
Its face has codimension one, so it
contributes a linear term to $d(v_B)$.

\subsection{The full and geometric algebras}
The passage from a square-zero derivation to higher Lie brackets
uses the dual form of the standard coderivation description of
homotopy algebras; see \cite[Chapter~10]{LodayVallette}.
The secondary-polytope construction and its geometric subalgebra
are those of \cite{KKS}, under the factorization hypothesis stated
above.

The polynomial differential just constructed
determines an $L_\infty$-algebra on
\begin{equation}\label{eq:full-linf}
\widetilde{\mathfrak g}_A
=W_A^\vee[-1]
=\bigoplus_{\substack{B\subseteq A\\\dim Q_B=d}}
\operatorname{or}(\Sigma(B))^\vee[-w(B)-1].
\end{equation}
This is the Chevalley--Eilenberg convention
\[
\bigl(\operatorname{Sym}
(\widetilde{\mathfrak g}_A^\vee[-1]),d_{\rm CE}\bigr)
=\bigl(\operatorname{Sym}(W_A),d\bigr).
\]
In particular, the line indexed by $B$ has
degree $w(B)+1=|B|-d$.

We specify the duality convention used below.
Let $U=W_A^\vee=\widetilde{\mathfrak g}_A[1]$,
and equip the symmetric coalgebra $S^c(U)$
with its unshuffle coproduct. Pair it with
$\operatorname{Sym}(W_A)$ by unnormalised graded
symmetrization. Thus, for distinct even basis
vectors and their duals,
$\langle v_1\cdots v_r,z_1\cdots z_r\rangle=1$.
The degree-one coderivation $\mathcal Q$ is
determined by
\begin{equation}\label{eq:duality-sign}
\langle d\alpha,c\rangle
=(-1)^{|\alpha|+1}\langle\alpha,\mathcal Qc\rangle.
\end{equation}
Its Taylor maps are
$\mathcal Q_r:S^r(U)\to U$.
There is no constant Taylor map, and
$\mathcal Q^2=0$ is equivalent to $d^2=0$.
This uses duality separately in each symmetric
power, so no completion is needed.

For the degree-$-1$ shift map
$s:\widetilde{\mathfrak g}_A\to U$, our skew
brackets are
\[
\ell_r=s^{-1}\mathcal Q_r s^{\otimes r},
\qquad |\ell_r|=2-r.
\]
All tensor products of graded maps include
their Koszul signs. This formula and
\eqref{eq:duality-sign} fix the desuspension
convention. The identities for the $\ell_r$
are equivalently the single equation
$\mathcal Q^2=0$.

The geometric subspace is
\begin{equation}\label{eq:geometric-linf}
\gA=
\bigoplus_{\substack{Q\text{ full-dimensional}\\
                    \operatorname{Vert}(Q)\subseteq A}}
\operatorname{or}(\Sigma(A\cap Q))^\vee
[-|A\cap Q|+d].
\end{equation}
Each summand corresponds to exactly one
geometric marked polytope.

\begin{proposition}\label{prop:geometric-subalgebra}
The subspace $\gA\subset\widetilde{\mathfrak g}_A$
is an $L_\infty$-subalgebra and $\ell_1|_{\gA}=0$.
\end{proposition}

\begin{proof}
A nonzero matrix coefficient of $\ell_r$
comes from a coarse subdivision
$\{(Q_i,B_i)\}_{i=1}^r$ of an output
$(Q_B,B)$. If every input is geometric,
then the covering $Q_B=\bigcup_i Q_i$ gives
\[
A\cap Q_B=\bigcup_i(A\cap Q_i)
=\bigcup_iB_i\subseteq B\subseteq A\cap Q_B.
\]
Thus the output is geometric as well.
For $r=1$, input and output have the same
underlying polytope, and a proper one-cell
subdivision has strictly smaller input
marking. A geometric input already contains
every ambient mark in that polytope, so
no such term is possible.
\end{proof}

The direction of the unary map is worth
making explicit: in the full $L_\infty$-algebra
it adds a missing mark to the input.
Its dual, the linear part of
\eqref{eq:cellular-derivation}, removes that
mark. Passing to \eqref{eq:geometric-linf}
eliminates these unary operations.

Put $B_A=\gA[1-d]$. The line associated to a
geometric polytope $Q$ now has degree
\begin{equation}\label{eq:shifted-degree}
|b_Q|=d+w(Q)=|A\cap Q|-1,
\qquad w(Q)=|A\cap Q|-d-1.
\end{equation}
The degree-$-d$ shift map
$s_d:B_A\to\gA[1]$ transports the Taylor maps
to
\begin{equation}\label{eq:tau-suspension}
\tau_r=s_d^{-1}\mathcal Q_r s_d^{\otimes r},
\qquad |\tau_r|=1-d(r-1).
\end{equation}
Explicitly, the tensor map in this formula
has sign
\[
s_d^{\otimes r}(b_1\otimes\cdots\otimes b_r)
=(-1)^{\,d\sum_i(r-i)|b_i|}
s_d b_1\otimes\cdots\otimes s_d b_r.
\]
Since the $\mathcal Q_r$ are graded symmetric,
\begin{equation}\label{eq:tau-symmetry}
\tau_r(\ldots,b_{i+1},b_i,\ldots)
=(-1)^{d+|b_i||b_{i+1}|}
\tau_r(\ldots,b_i,b_{i+1},\ldots).
\end{equation}
Equivalently, in the endomorphism operad,
$\tau_r\cdot\sigma=(\sgn\sigma)^d\tau_r$.
We use this suspension convention when
normalizing the fundamental-chain morphism.

\begin{lemma}\label{lem:weight-growth}
Every nonzero primary operation $\tau_r$,
$r\geq2$, satisfies
\[
w_{\rm out}=1+\sum_iw_{{\rm in},i}.
\]
More generally, a nonzero evaluation of
a tree with $k$ primary-operation vertices
has output weight equal to the sum of
its input weights plus $k$.
\end{lemma}

\begin{proof}
The first assertion is \eqref{eq:facet-weight};
it also follows from
\eqref{eq:shifted-degree} and
\eqref{eq:tau-suspension}.
Apply this identity successively at the
vertices of the tree to obtain the second.
\end{proof}

All weights lie between zero and $|A|-d-1$.
Thus every composition of three primary
operations vanishes when $|A|\leq d+3$.
This statement allows nonzero compositions
of two operations; it supplies the finite
weight bound used in the later discussion.

\subsection{Circuit operations}
We recall the two triangulations of an affine circuit, as in
\cite{GKZ,DLRS}, and record the resulting operations with our
orientation and suspension conventions.

Suppose now that
$A=\{a_1,\ldots,a_{d+2}\}$ is a full-dimensional
affine circuit. Choose its unique affine
relation up to scale,
\[
\sum_i\lambda_i a_i=0,\qquad
\sum_i\lambda_i=0,\qquad \lambda_i\neq0,
\]
and put $I_\pm=\{i:\pm\lambda_i>0\}$.
Write $B_i=A\setminus\{a_i\}$ and
$\Delta_i=\Conv(B_i)$.

\begin{proposition}\label{prop:circuit-triangulations}
The circuit has exactly two regular
triangulations,
\[
\mathcal T_+=\{(\Delta_i,B_i):i\in I_+\},
\qquad
\mathcal T_-=\{(\Delta_i,B_i):i\in I_-\}.
\]
They are the endpoints of the secondary
interval $\Sigma(A)$.
\end{proposition}

\begin{proof}
For heights $h_i$, set
$\delta(h)=\sum_i\lambda_i h_i$.
Let $h^{(i)}$ be the affine function agreeing
with the heights on $B_i$. The affine
relation gives
\[
h_i-h^{(i)}(a_i)=\frac{\delta(h)}{\lambda_i}.
\]
For $\delta(h)\neq0$, the lifted simplex
on $B_i$ is a lower facet precisely when
this quantity is positive. Its projections
are therefore the cells of $\mathcal T_+$
or $\mathcal T_-$ according to the sign of
$\delta(h)$. If $\delta(h)=0$, all lifted
points lie on the graph of an affine
function, giving the trivial subdivision.
Since $\dim\Sigma(A)=1$, the two
triangulations are its endpoints.
\end{proof}

Orient $\Sigma(A)$ from $\mathcal T_-$ to
$\mathcal T_+$. The simplex secondary
polytopes are oriented points, so the
full cellular differential is
\begin{equation}\label{eq:circuit-differential}
d(v_A)=\prod_{i\in I_+}v_{B_i}
       -\prod_{i\in I_-}v_{B_i},
\qquad d(v_{B_i})=0.
\end{equation}
If one sign set is a singleton, one
monomial here is linear.

Let $t_i$ be the shifted dual generator
of an empty $\Delta_i$, and let $q$
be the shifted dual of $v_A$.
The geometric space is
\[
B_A^d=\bigoplus_{\Delta_i\ {\rm empty}}\R t_i,
\qquad
B_A^{d+1}=\R q,
\qquad B_A^j=0\quad(j\neq d,d+1).
\]
If both $|I_+|$ and $|I_-|$ are at least
two, every $\Delta_i$ is empty.
Indeed, containment $a_i\in\Delta_i$
would give an affine relation whose
Radon partition has singleton side
$\{i\}$, contradicting uniqueness of
the circuit relation.
If $I_-=\{j\}$, then $a_j$ is interior
to $\Delta_j$, so $\Delta_j$ is not
geometric with its vertex marking.
Every other deletion simplex is empty.
The case $|I_+|=1$ is analogous.

In the dual bases specified above, put
$\epsilon_{d,r}=(-1)^{d\binom r2}$.
Equations \eqref{eq:duality-sign} and
\eqref{eq:tau-suspension} give
\begin{equation}\label{eq:circuit-operations}
\begin{aligned}
\tau_{|I_+|}\bigl((t_i)_{i\in I_+}\bigr)
   &=\epsilon_{d,|I_+|}\,q,
       && |I_+|\geq2,\\
\tau_{|I_-|}\bigl((t_i)_{i\in I_-}\bigr)
   &=-\epsilon_{d,|I_-|}\,q,
       && |I_-|\geq2.
\end{aligned}
\end{equation}
To see the sign directly, the simplex
generators in $U$ have degree zero,
while $|v_A|=-1$. Thus
\eqref{eq:duality-sign} introduces no
sign on \eqref{eq:circuit-differential}.
The remaining tensor-suspension sign
is $(-1)^{d^2\binom r2}
=\epsilon_{d,r}$.
All permutations of a displayed tuple
have the same value, by
\eqref{eq:tau-symmetry} and $|t_i|=d$.
All other basis tuples have value zero.
When the two arities agree, the two
lines specify different entries of
the same operation.

The simplex generators appearing in
the two triangulations are disjoint.
Consequently, changing the sign of one
generator in each nonunary tuple, when
needed, gives the simpler normalization
with coefficients $+1$ and $-1$.
We will use this normalization in the
examples.

\begin{corollary}\label{cor:circuit-central}
For a circuit, let $Z_A=\R q$.
Every $\tau_r$, $r\geq2$, has image
in $Z_A$ and vanishes whenever an
input belongs to $Z_A$. In particular,
\[
\tau_k\circ_i\tau_l=0
\qquad(k,l\geq2).
\]
\end{corollary}

\begin{proof}
Proper full-dimensional subconfigurations
of a circuit are simplices, whose
secondary polytopes have no facets.
Thus \eqref{eq:circuit-operations}
lists every nonzero higher operation.
Each has output $q$ and only simplex
inputs. Equivalently, two primary
operations would require weight at
least two, whereas a circuit has
weights zero and one.
\end{proof}

This is the property needed for the
configuration-space construction of
Section~\ref{sec:configuration-space-action}. The circuit $L_\infty$-algebra
has now been specified, including the
effect of geometric markings, the
degree shift, and the suspension signs.

\section{Operad models and the algebraic enhancement}
\label{sec:operad-models}
\label{sec:abstract-enhancement}
We use two dg operads modeling little $d$-disks. Geometric chains
provide a direct formulation of the enhancement problem, while the
homotopy Poisson operad gives its algebraic solution.
For the little cubes viewpoint see \cite{May}; for dg operads,
their homotopy algebras, and the bar--cobar formalism see
\cite{LodayVallette}.

\subsection{Configuration spaces and geometric chains}
For $r\geq2$, let
\[
C_d(r)=
\operatorname{Conf}_r(\R^d)/
(\R^d\rtimes\R_{>0}),
\]
where translations and positive dilations act
simultaneously on the ordered configuration.
Each orbit has a unique representative satisfying
\[
\sum_{i=1}^r x_i=0,\qquad
\sum_{i=1}^r\lVert x_i\rVert^2=1.
\]
Thus $C_d(r)$ is an open submanifold of a sphere
of dimension $D_r=d(r-1)-1$.
For $d\geq2$ it is connected.

The compactification originates in \cite{FM}. We use its real
manifold-with-corners version, described by Sinha \cite{Sinha}
and in the operadic setting of \cite[Section~5]{LV}.
The Fulton--MacPherson compactification
$M_r=\FM_d(r)$ records the relative configurations
of points in colliding clusters. It is a compact
smooth semialgebraic manifold with corners,
with interior $C_d(r)$.
Set $M_1=\{*\}$ and omit arity zero.
For $k,l\geq2$, the insertion map
\begin{equation}\label{eq:fm-insertion}
\gamma_i:M_k\times M_l\longrightarrow M_{k+l-1}
\end{equation}
identifies its source with the closed face where
the $l$ inserted points form an infinitesimal
cluster. In particular, its image lies in the
boundary and
\[
D_k+D_l=D_{k+l-1}-1.
\]
Higher boundary strata correspond to nested
clusters, or equivalently to rooted trees.
These spaces and insertion maps form a model
of the little $d$-disks operad
\cite[Section~5]{LV}.

Permutations act on the orientation of $M_r$
by the character
\begin{equation}\label{eq:fm-orientation-character}
\chi_d(\sigma)=(\sgn\sigma)^d.
\end{equation}
Indeed, permuting the $r$ blocks of dimension
$d$ in $(\R^d)^r$ has this determinant.
The diagonal translation subspace is fixed
pointwise, and the radial direction is
preserved. The same character therefore acts
on the normalized configuration space and its
compactification. Reversing the chosen
orientation in a fixed arity does not change
this character.

Piecewise algebraic chains and forms on compactified configuration
spaces were used by Kontsevich--Soibelman \cite[Appendix]{KSDeligne}.
We use the geometric chains of \cite[Section~3]{HLTV},
with real coefficients. Concretely, a chain is
a finite linear combination of integration
currents of oriented semialgebraic submanifolds
with compact closure, with equality as currents.
Their boundary, pushforward, and cross product
give a chain complex naturally equivalent to
singular chains. In particular,
\[
\partial(c\times c')
=\partial c\times c'+(-1)^{\dim c}c\times\partial c'.
\]
These chains include the fundamental chains
of compact oriented manifolds with corners
and their faces.

Place a $p$-chain in cohomological degree $-p$,
and define
\[
\En(r)=C^{\mathrm{SA}}_{-*}(M_r;\R),
\qquad d_{\En}=\partial.
\]
The composition is
$c\circ_i c'=(\gamma_i)_*(c\times c')$.
The chain of the unary point is the operadic
identity. This defines the positive-arity
dg operad used for the geometric construction. Omitting arity zero
means that the algebra has no specified multiplicative unit; the
unary operadic identity still acts as the identity map. This is
the nonunital convention in the disk-algebra comparison of
\cite[Proposition~5.4.5.15]{HA}.

Smooth forms on $M_r$ can be evaluated on these
currents. We will use
\begin{equation}\label{eq:chain-integration}
\int_{\partial c}\omega=\int_c d\omega,
\qquad
\int_{f_*c}\omega=\int_c f^*\omega,
\end{equation}
for the smooth maps occurring here.
The auxiliary forms below are ordinary smooth
forms; their integration does not require
them to be semialgebraic. This follows directly
from the current description; see also
\cite[Section~5.1]{HLTV}.

\subsection{The fundamental-chain morphism}
Let $\mathsf L_{\infty,d}$ denote the suspended
$L_\infty$ operad governing the operations
on $V$ obtained from an $L_\infty$-structure on
$V[d-1]$. Its free graded-operad generators
$\mu_r$, $r\geq2$, have degree
\[
|\mu_r|=1-d(r-1)=-D_r
\]
and symmetric-group character $\chi_d$.
The differential of a generator is a signed
sum of two-vertex trees. Its algebra
relations are the square-zero coderivation
identities, transported by
\eqref{eq:tau-suspension}. The differential
of the underlying complex is part of the
dg-algebra structure, rather than an
additional free unary generator.

We use the standard fundamental-chain morphism
\begin{equation}\label{eq:fundamental-chain-map}
\iota_d:\mathsf L_{\infty,d}\longrightarrow\En,
\qquad
\mu_r\longmapsto[M_r].
\end{equation}
For this construction see
\cite[Section~2.2, equation~(1)]{CW}.
We transport its generator convention to
\eqref{eq:tau-suspension} and choose the
orientation of each $M_r$ accordingly.
This fixes the aritywise signs in
\eqref{eq:fundamental-chain-map}.
All subsequent integrals use these orientations.

The chain identity underlying this morphism
is
\begin{equation}\label{eq:fundamental-boundary}
\partial[M_r]
=\sum_{\substack{T\text{ a rooted tree}\\
                 \text{with two internal vertices}}}
\varepsilon(T)(\gamma_T)_*
\bigl([M_{k(T)}]\times[M_{l(T)}]\bigr),
\end{equation}
where the leaves are labelled and each boundary
face is counted once. The sign
$\varepsilon(T)$ compares its product
orientation with the boundary orientation.
Permuting labels and exchanging oriented
factors supplies the suspension and Koszul
signs in the differential of
$\mathsf L_{\infty,d}$. Thus
\eqref{eq:fundamental-boundary} expresses
the chain-map property of $\iota_d$.

For $r\geq3$, $[M_r]$ is a fundamental
\emph{chain with boundary}, representing
the relative fundamental class of
$(M_r,\partial M_r)$. It is not asserted
to be an absolute cycle. For $r=2$,
$M_2=S^{d-1}$ and its fundamental chain
is a cycle.

\subsection{The homotopy Poisson retraction}
Let $\mathsf{Pois}_d$ denote the positive-arity Poisson operad
with a degree-zero commutative product and a Lie bracket of
cohomological degree $1-d$. We use the cobar models of
\cite[Section~2]{TW}, with homological degrees reversed. In the
notation of that reference, $\{d\}$ is an operadic suspension,
$(-)^*$ is aritywise linear duality, and $\Omega$ is the reduced
operadic cobar construction. General references for these
constructions and Koszul resolutions are
\cite[Chapters~6--7]{LodayVallette}. With this convention the models are
\[
 P_d^\infty=\Omega\bigl(\mathsf{Pois}_d^*\{d\}\bigr),
 \qquad
 L_d^\infty=\Omega\bigl(\mathsf{Com}^*\{d\}\bigr).
\]
Here $L_d^\infty$ is $\mathsf L_{\infty,d}$ in the suspension
convention of \eqref{eq:tau-suspension}. The primary Lie generators
have degree $1-d(r-1)$. The Koszul resolution
$P_d^\infty\to\mathsf{Pois}_d$ is a quasi-isomorphism. For
$d\geq2$ the homology operad of little $d$-disks is
$\mathsf{Pois}_d$, and real operadic formality connects it to the
chain operad by quasi-isomorphisms \cite[Theorem~1.1]{LV}.
The homology calculation is due to Cohen \cite[Part~III]{CLM};
for the formality construction see also \cite{Kontsevich}.
Restricting to positive arities gives the model used here.

\begin{lemma}\label{lem:poisson-retraction}
The shifted Lie inclusion $i:L_d^\infty\to P_d^\infty$ admits
a dg-operad retraction $p:P_d^\infty\to L_d^\infty$ which sends
the binary product generator to zero.
\end{lemma}

\begin{proof}
Consider the operad maps
\[
 \mathsf{Com}\xrightarrow{u}\mathsf{Pois}_d
     \xrightarrow{v}\mathsf{Com},\qquad vu=1.
\]
The map $u$ includes the commutative product. The map $v$ preserves
that product and kills the bracket, so it respects the Jacobi and
Leibniz relations. These operads are finite-dimensional in each
arity. Dualizing, suspending and applying $\Omega$ gives
\[
 i=\Omega(v^*\{d\}),\qquad p=\Omega(u^*\{d\}),\qquad pi=1.
\]
The generators dual to the pure commutative monomials span the
shifted homotopy Lie generators in the cobar model \cite[Section~2]{TW}.
In arity two the other generator, dual to the bracket, represents
the degree-zero product. Its restriction along $u$ vanishes.
\end{proof}

\begin{proof}[Proof of Proposition~\ref{prop:intro-algebraic}]
Write $\tau:L_d^\infty\to\End_{B_A}$ for the action obtained in
Section~\ref{sec:secondary-polytopes}. Define
\[
     \rho_A^{\mathrm{alg}}=\tau p:
       P_d^\infty\longrightarrow\End_{B_A}.
\]
Then $\rho_A^{\mathrm{alg}}i=\tau$ as dg-operad maps. This recovers
each of the chosen higher operations on $B_A$ itself. The binary
product vanishes by Lemma~\ref{lem:poisson-retraction}. The proof
uses no upper bound on the weights or on $|A|$.
\end{proof}

The equality just proved uses $i$ in the cobar model. The geometric
construction uses the fundamental-chain map $\iota_d:L_d^\infty\to\En$
of \cite[Section~2.2]{CW}. The two inclusions agree after comparison
in the homotopy category, with the following choice of normalization.
Choose a formality isomorphism from $\En$ to $\mathsf{Pois}_d$ in
that category which sends $[M_2]=[S^{d-1}]$ to the Poisson bracket.
This normalization is possible by rescaling the bracket. Since
$L_d^\infty$ is cofibrant and $\mathsf{Pois}_d$ is fibrant,
the composite with $\iota_d$ has a dg-operad representative
$L_d^\infty\to\mathsf{Pois}_d$ \cite{BM}. Its binary generator
maps to the bracket. For $r>2$ its generator has degree
\[
    1-d(r-1)=(1-d)(r-1)+(2-r)<(1-d)(r-1).
\]
The last quantity is the lowest degree in $\mathsf{Pois}_d(r)$,
since an $r$-ary monomial contains at most $r-1$ brackets.
All the higher generators must therefore map to zero. This
determines the representative on every generator and identifies
it with $L_d^\infty\xrightarrow{i}P_d^\infty\to\mathsf{Pois}_d$.

The following sections realize the operations on these chains by
integration. The insertion faces determine the boundary values of
the forms used in the construction.

\section{Geometric enhancements for circuits}
\label{sec:configuration-space-action}
We first treat operations whose images are annihilated by all further
operations. The circuit algebra has this property. In this case a
normalized form in each arity defines the full action.

\subsection{Normalized forms}
\begin{lemma}\label{lem:normalized-forms}
For every $r\geq2$, there is a smooth
$D_r$-form $\eta_r$, compactly supported
in the interior of $M_r$, such that
\begin{equation}\label{eq:form-normalization}
\int_{M_r}\eta_r=1,\qquad
R_\sigma^*\eta_r=\chi_d(\sigma)\eta_r
\quad(\sigma\in S_r),
\end{equation}
where $R_\sigma$ is the relabelling map.
\end{lemma}

\begin{proof}
Choose a top-degree form $\eta_r^0$
supported in an oriented coordinate ball
in the interior, with integral one.
Set
\[
\eta_r=\frac1{r!}
\sum_{\sigma\in S_r}
\chi_d(\sigma)R_\sigma^*\eta_r^0.
\]
The finite average has compact support
in the interior and transforms by
$\chi_d$. By
\eqref{eq:fm-orientation-character},
each summand, including its coefficient,
has integral one. The average therefore
has integral one. It is closed because
it has top degree.
\end{proof}

We fix these forms once for each $d$ and
$r$. For $r=2$, one may take a normalized
volume form on $S^{d-1}$, with the orientation
fixed above. The support condition is then
automatic because the sphere has no boundary.

\subsection{The extension lemma}
\begin{theorem}\label{thm:central-extension}
Let $V$ be a graded real vector space with
zero differential, and let $Z\subseteq V$
be a graded subspace. Suppose that, for
each $r\geq2$, a homogeneous multilinear map
\[
\tau_r:V^{\otimes r}\longrightarrow V
\]
satisfies the following conditions:
\begin{enumerate}
\item $|\tau_r|=-D_r$ and
      $\tau_r\cdot\sigma=\chi_d(\sigma)\tau_r$;
\item $\operatorname{im}\tau_r\subseteq Z$;
\item $\tau_r(v_1,\ldots,v_r)=0$ whenever
      some $v_i$ belongs to $Z$.
\end{enumerate}
Then the formula
\begin{equation}\label{eq:central-chain-action}
\rho_r(c)=
\begin{cases}
\displaystyle\left(\int_c\eta_r\right)\tau_r,
     &\dim c=D_r,\\[5pt]
0,   &\dim c\neq D_r
\end{cases}
\qquad(r\geq2),
\end{equation}
together with $\rho_1([*])=\operatorname{id}_V$,
defines a morphism of dg operads
\[
\rho:\En\longrightarrow\End_V.
\]
Its restriction along $\iota_d$ sends
$\mu_r$ to $\tau_r$.
\end{theorem}

\begin{proof}
Integration is well-defined on geometric
chains because they are identified as
currents. We verify the dg-operad identities.

\emph{Degrees and differential.}
A $D_r$-chain has degree $-D_r$, the
degree of $\tau_r$. Hence $\rho_r$ has
degree zero. The differential of
$\End_V$ is zero. The only possibly
nonzero value of $\rho_r(\partial c)$
would require $\dim c=D_r+1$, and
\eqref{eq:chain-integration} gives
\[
\left(\int_{\partial c}\eta_r\right)\tau_r
=\left(\int_c d\eta_r\right)\tau_r=0.
\]
In the geometric chain model such chains
already vanish by dimension. All other
degrees give zero directly. Thus
$\rho_r\partial=0$.

\emph{Equivariance.}
For a $D_r$-chain, naturality of integration
and \eqref{eq:form-normalization} give
\[
\rho_r((R_\sigma)_*c)
=\left(\int_cR_\sigma^*\eta_r\right)\tau_r
=\chi_d(\sigma)\rho_r(c)
=\rho_r(c)\cdot\sigma.
\]
Both sides are zero in the remaining
degrees.

\emph{Composition.}
Let $k,l\geq2$ and $m=k+l-1$.
The map $\gamma_i$ takes values in
$\partial M_m$. Since $\eta_m$ is
supported in the interior,
$\gamma_i^*\eta_m=0$, and therefore
\[
\rho_m(c\circ_i c')=0
\]
for arbitrary homogeneous chains
$c\in\En(k)$ and $c'\in\En(l)$.
On the endomorphism side, a nonzero
$\rho_l(c')$ has image in $Z$, which
is annihilated in each input of
$\rho_k(c)$. Consequently,
\[
\rho_k(c)\circ_i\rho_l(c')=0.
\]
The two sides agree, including all
degrees in which one of the actions
is already zero. Compositions with
the unary point agree because it
acts as the identity.

This proves every dg-operad identity.
Finally,
\[
\rho_r([M_r])
=\left(\int_{M_r}\eta_r\right)\tau_r
=\tau_r,
\]
which proves the restriction statement.
The prescribed maps do define a
$\mathsf L_{\infty,d}$-algebra: every
two-vertex composite is zero by
conditions (2) and (3), so the
generator differentials act by zero.
\end{proof}

\begin{remark}\label{rem:zero-product}
Since $D_2=d-1>0$, every zero-dimensional
binary chain acts by zero in this
construction. Thus its ordinary product
vanishes. The unary identity remains
the identity of $V$. The theorem uses
the positive-arity operad, as specified
in the introduction.
\end{remark}

For the fixed family of forms, the construction
is compatible with strict morphisms. If
$f:V\to V'$ is a degree-zero linear map with
$f\tau_r=\tau'_r f^{\otimes r}$ for all $r$,
then \eqref{eq:central-chain-action} gives
$f\rho_r(c)=\rho'_r(c)f^{\otimes r}$ for
every chain $c$.

\subsection{Changing the auxiliary forms}
The existence statement involves choices of
forms. Within this construction, their effect
can be controlled while fixing all the
fundamental-chain operations.

\begin{proposition}\label{prop:form-homotopy}
Let $\eta_r^0$ and $\eta_r^1$ be two families
satisfying Lemma~\ref{lem:normalized-forms},
and let $\rho^0,\rho^1$ be the resulting
actions on the same $(V,\{\tau_r\})$.
There is a polynomial homotopy of
dg-operad maps
\[
\mathcal H:\En\longrightarrow
\End_V\otimes\R[t,dt]
\]
with endpoints $\rho^0,\rho^1$ and
\[
(\mathcal H\circ\iota_d)(\mu_r)=\tau_r\otimes1.
\]
Here $|t|=0$, $|dt|=1$, and $d(t)=dt$.
\end{proposition}

\begin{proof}
For each arity, the compactly supported
top form $\eta_r^1-\eta_r^0$ has integral
zero on the connected oriented manifold
$C_d(r)$. Integration identifies
$H_c^{D_r}(C_d(r);\R)$ with $\R$.
Consequently there is a compactly supported
smooth $(D_r-1)$-form $\beta_r$ on the
interior satisfying
\[
d\beta_r=\eta_r^1-\eta_r^0.
\]
Averaging $\beta_r$ with the character
$\chi_d$ makes it equivariant without
changing this equation. Extend it by
zero to $M_r$ and put
$\eta_r(t)=(1-t)\eta_r^0+t\eta_r^1$.

For a homogeneous $p$-chain $c$, define
\begin{equation}\label{eq:form-homotopy}
\mathcal H_r(c)=
\begin{cases}
\displaystyle
\tau_r\otimes\left(\int_c\eta_r(t)\right),
   &p=D_r,\\[6pt]
\displaystyle
(-1)^{D_r}\tau_r\otimes
        \left(\int_c\beta_r\right)dt,
   &p=D_r-1,\\[6pt]
0, &\text{otherwise}.
\end{cases}
\end{equation}
The unary point maps to
$\operatorname{id}_V\otimes1$.
Each line has total degree $-p$.

For $p=D_r$, the tensor-product
differential and Stokes' formula give
\[
\begin{aligned}
d\,\mathcal H_r(c)
&=(-1)^{D_r}\tau_r\otimes
   \left(\int_c(\eta_r^1-\eta_r^0)\right)dt\\
&=(-1)^{D_r}\tau_r\otimes
   \left(\int_{\partial c}\beta_r\right)dt
 =\mathcal H_r(\partial c).
\end{aligned}
\]
For $p=D_r-1$ both sides are zero, and
all remaining degrees follow as in
Theorem~\ref{thm:central-extension}.
Equivariance follows from that of
$\eta_r^0,\eta_r^1,\beta_r$.

Both $\eta_r(t)$ and $\beta_r$ vanish
near the boundary. Hence $\mathcal H$
vanishes on every composite of two
nonunary chains. The corresponding
composition in $\End_V\otimes\R[t,dt]$
is zero because its endomorphism
factor contains a composite of two
$\tau$'s. This proves the operadic
composition identity.

Evaluation at $t=0$ or $t=1$ sets
$dt=0$ and gives the required endpoints.
Finally, both top forms have integral
one, so
\[
\mathcal H_r([M_r])=\tau_r\otimes1.
\]
Since the $\mu_r$ generate the
underlying graded operad, the homotopy
is constant on its prescribed
$\mathsf L_{\infty,d}$-action.
\end{proof}

The homotopy compares the actions obtained from
\eqref{eq:central-chain-action}, with the prescribed higher Lie
operations fixed.

\subsection{Application to secondary-polytope circuits}
\begin{corollary}[Circuit enhancement]\label{cor:circuit-enhancement}
For a full-dimensional affine circuit $A\subset\R^d$, formula
\eqref{eq:central-chain-action} defines an $\En$-action on $B_A$
recovering its primary operations on fundamental chains.
\end{corollary}
\begin{proof}
By Proposition~\ref{prop:geometric-subalgebra},
the geometric circuit algebra has zero
differential. Its shifted operations have
the required degrees and symmetry by
\eqref{eq:tau-suspension} and
\eqref{eq:tau-symmetry}.
Corollary~\ref{cor:circuit-central}
supplies the subspace $Z_A=\R q$ and
the image and vanishing conditions of
Theorem~\ref{thm:central-extension}.
That theorem therefore defines
\[
\rho_A:\En\longrightarrow\End_{B_A}
\]
on the original graded space, with
$(\rho_A\circ\iota_d)(\mu_r)=\tau_r$.
Shifting back by $[d-1]$ recovers
the geometric secondary-polytope
$L_\infty$-algebra in the convention
of Section~\ref{sec:secondary-polytopes}.

The description of the graded space
was established in the circuit
calculation. Its degrees are $d$
and $d+1$, so any degree-zero
binary operation must vanish:
two inputs have total degree at
least $2d>d+1$. This also follows
from Remark~\ref{rem:zero-product}.
\end{proof}

\begin{corollary}\label{cor:circuit-homotopy}
For a fixed circuit algebra and fixed
orientation and suspension conventions,
the actions obtained by varying the
normalized forms are homotopic relative
to their shifted $L_\infty$-structure.
\end{corollary}

\begin{proof}
Apply Proposition~\ref{prop:form-homotopy}
with $V=B_A$ and $Z=Z_A$.
\end{proof}

\section{The weight filtration and extension obstructions}
\label{sec:weight-extension}

The convex pentagon requires operations on chains
of codimension one in a configuration space.
More generally, a correction of weight $w$
will be represented by a form of degree
$d(r-1)-w$ in arity $r$.
We describe the extension problem for these
forms. Reflection on the auxiliary configuration spaces allows
us to carry the construction through weight $2d-1$.

\subsection{Closed extensions from the boundary}
Let $X$ be a compact smooth manifold with
corners whose boundary hypersurfaces are embedded.
Write $\Omega^*(\partial X)$ for the complex
of collections of smooth forms on its boundary
hypersurfaces whose pullbacks agree on every
intersection. These are tangential boundary
values; no normal derivatives are prescribed.
Let
\[
\Omega^*(X,\partial X)
=\ker\bigl(\Omega^*(X)\longrightarrow
                 \Omega^*(\partial X)\bigr).
\]
Both complexes may be tensored with a
finite-dimensional real vector space $E$.

\begin{lemma}\label{lem:relative-form-extension}
For a closed boundary form
$K\in\Omega^p(\partial X;E)$, $p\geq0$,
the obstruction to a closed extension with
the same boundary values is
\[
\delta[K]\in H^{p+1}(X,\partial X;\R)\otimes E.
\]
If this class vanishes, the closed extensions,
modulo exact forms with a primitive vanishing
on the boundary, form a torsor for
$H^p(X,\partial X;\R)\otimes E$.
Equivariant boundary data for a finite group
admit an equivariant closed extension whenever
they admit a closed extension.
\end{lemma}

\begin{proof}
Restriction gives a short exact sequence
\[
0\longrightarrow\Omega^*(X,\partial X;E)
\longrightarrow\Omega^*(X;E)
\longrightarrow\Omega^*(\partial X;E)
\longrightarrow0.
\]
To see surjectivity, work in a corner chart
$[0,\infty)^a\times\R^b$.
For each differential-form component, extend
its compatible coefficient functions from the
coordinate faces by inclusion and exclusion.
Explicitly, for a component with differential
$dx_I$, let $L$ be the set of normal coordinate
indices not in $I$. Its prescribed coefficients
on the faces indexed by $L$ have common
restrictions $f_J$ on intersections indexed by
nonempty $J\subseteq L$. A local extension is
\[
\sum_{\varnothing\neq J\subseteq L}
(-1)^{|J|+1}\pi_J^*f_J,
\]
where $\pi_J$ sets the coordinates indexed
by $J$ equal to zero. Restricting to any
one of those faces cancels the extra terms
and recovers its prescribed coefficient.
Components containing a given
normal differential have zero pullback to that
face and impose no condition there.
A partition of unity gives a global extension.

The boundary complex computes the cohomology
of $\partial X$. Indeed, the local radial
contraction in a corner chart preserves all
faces and gives the Poincare lemma for
compatible collections; partitions of unity
give a fine resolution. The kernel therefore
computes relative cohomology.

Choose any extension $\widetilde K$.
Its differential vanishes on every face,
and the relative class $[d\widetilde K]$
is $\delta[K]$. This class is independent
of the extension. If it is zero, write
$d\widetilde K=d\xi$ with
$\xi\in\Omega^p(X,\partial X;E)$.
Then $\widetilde K-\xi$ is a closed
extension with the prescribed boundary values.
Conversely a closed extension makes the
class zero.

Two closed extensions differ by a closed
relative form. Taking their quotient by
exact relative forms gives the stated torsor.
Finally average a closed extension over the
finite group, using pullback on forms and
the inverse action on coefficients.
\end{proof}

If $X$ is connected and oriented of dimension
$D$, Poincare--Lefschetz duality identifies
the obstruction space, as a real vector space,
with
\begin{equation}\label{eq:extension-duality}
H^{p+1}(X,\partial X;\R)\otimes E
\cong H_{D-p-1}(X;\R)\otimes E.
\end{equation}
For $p=D-1$ the obstruction vanishes precisely
when
\begin{equation}\label{eq:top-boundary-integral}
\int_{\partial X}K=0.
\end{equation}
This follows either from duality or by
integrating $d\widetilde K$ and applying
Stokes' theorem. When $\partial X$ is empty,
the prescribed boundary value is zero and
the zero form is an available extension.

\subsection{Forms graded by weight}
Let
\[
V=\bigoplus_{a=0}^W V_{(a)}
\]
be a finite-dimensional graded real vector
space with zero differential and an additional
nonnegative weight grading.
An endomorphism has weight $w$ if it takes
inputs of weights $a_1,\ldots,a_r$ to weight
$a_1+\cdots+a_r+w$.
Maps of weight greater than $W$ vanish.

Suppose a fixed
$\mathsf L_{\infty,d}$-algebra structure on
$V$ has primary operations $\tau_r$ of
weight one. Set
\[
p_{r,w}=d(r-1)-w=D_r+1-w,\qquad
E_{r,w}=\End^{-p_{r,w}}_{(w)}(V)(r).
\]
We seek closed forms
\[
\alpha_{r,w}\in\Omega^{p_{r,w}}(M_r;E_{r,w}),
\qquad r\geq2,\quad 1\leq w\leq W.
\]
Forms of negative degree are understood to
be zero. At weight one put
\begin{equation}\label{eq:primary-weight-form}
\alpha_{r,1}=\eta_r\,\tau_r,
\end{equation}
with the normalization of
Lemma~\ref{lem:normalized-forms}.
At every weight $w\geq2$ we choose
$\alpha_{2,w}=0$; the binary configuration space
has no boundary.

We fix the insertion convention for
endomorphism-valued forms.
For a scalar $p$-form $\alpha$ and a
degree-$-p$ endomorphism $f$, and similarly
for $\beta$ and $g$, write
\[
(\alpha f)\star_i(\beta g)
=\bigl(\operatorname{pr}_1^*\alpha
       \wedge\operatorname{pr}_2^*\beta\bigr)
       (f\circ_i g).
\]
The factors are ordered outer, then inner.
The coefficient uses composition in the
graded endomorphism operad, including its
insertion sign. Thus integration over
a product chain gives the composition of
the two integrated endomorphisms.
This formula defines $\star_i$; no additional
tensor-product sign is implicit.

For a face
$\gamma_i:M_k\times M_l\longrightarrow M_r$,
$r=k+l-1$, the required boundary equation is
\begin{equation}\label{eq:weight-boundary-equation}
\gamma_i^* \alpha_{r,w}
=\sum_{\substack{a+b=w\\a,b\geq1}}
       \alpha_{k,a}\star_i \alpha_{l,b}.
\end{equation}
Relabellings are applied simultaneously
to the face and to the endomorphism.
The equation has the correct form degree,
since
\[
p_{k,a}+p_{l,b}=p_{r,a+b}.
\]
At weight one its right-hand side is zero,
as is the boundary restriction of
\eqref{eq:primary-weight-form}.

\begin{proposition}\label{prop:weight-obstruction}
Suppose the closed equivariant forms have
been constructed in all weights less than
$w$, satisfying
\eqref{eq:weight-boundary-equation}.
Its right-hand sides define compatible closed
boundary data $K_{r,w}$. For $p_{r,w}\geq0$
the obstruction to extending these data is
\begin{equation}\label{eq:weight-obstruction}
o_{r,w}=\delta[K_{r,w}]
\in H^{D_r+2-w}(M_r,\partial M_r;\R)
       \otimes E_{r,w}.
\end{equation}
After fixing the orientation, the space in
\eqref{eq:weight-obstruction} is isomorphic to
\[
H_{w-2}(M_r;\R)\otimes E_{r,w}.
\]
The class depends on the lower-weight forms
already chosen.
\end{proposition}

\begin{proof}
Every summand in the boundary prescription
is closed. At a codimension-two corner,
restriction in either order expands into
the same sum of three-vertex expressions,
with vertex weights adding to $w$.
For nested insertions this follows from
associativity in the endomorphism operad.
For disjoint insertions, exchanging factors
of form degrees $p$ and $q$ contributes
$(-1)^{pq}$. The interchange of the
degree-$-p$ and degree-$-q$ endomorphisms
contributes the same sign. Their product
is $+1$, so the two restrictions agree.
The same argument on further intersections
gives compatible boundary data.
Relabelling commutes with the construction,
so the data are equivariant.

Apply Lemma~\ref{lem:relative-form-extension}
with $p=p_{r,w}$. Substituting this degree
in \eqref{eq:extension-duality} gives
$D_r-p_{r,w}-1=w-2$.
If $p_{r,w}<0$, all boundary expressions
also have negative total form degree and
are zero, so there is no extension problem.
\end{proof}

\subsection{Reflection on configuration spaces}
Fix a linear hyperplane reflection $s$ of the auxiliary $\R^d$.
It acts on $M_r$ by reflecting all $r$ points, commutes with
relabelling, and commutes with every insertion when it acts on
both factors. There is no symmetry assumption on the marked set
$A$. We let $s$ act trivially on the coefficient spaces $E_{r,w}$.
Its orientation character on $M_r$ is
\[
     \epsilon_r=(-1)^{r-1}.
\]
Indeed, centered configurations form an open set in the unit
sphere of $(\R^d)^r/\R^d_{\mathrm{diag}}$. The determinant on the
quotient is $(\det s)^r/\det s$. The outward radial line is
preserved, so this is also the degree on the oriented interior and
on the relative fundamental class of the compactification.
Since $\epsilon_{k+l-1}=\epsilon_k\epsilon_l$, the condition
\begin{equation}\label{eq:reflection-type}
         s^*\alpha_{r,w}=\epsilon_r \alpha_{r,w}
\end{equation}
is preserved by the boundary equation \eqref{eq:weight-boundary-equation}.
For the primary forms we replace $\eta_r$ by
$(\eta_r+\epsilon_r s^*\eta_r)/2$. This preserves its integral,
interior support and permutation character. At a later stage,
whenever a closed extension exists, the same averaging gives an
extension of type \eqref{eq:reflection-type} with its prescribed
boundary values. The permutation and reflection averages commute.

\begin{lemma}\label{lem:reflection-obstruction}
Assume the lower-weight forms satisfy \eqref{eq:reflection-type}.
Under Poincar\'e--Lefschetz duality the obstruction $o_{r,w}$
belongs to
\[
      H_{w-2}(M_r;\R)^{s=+1}\otimes E_{r,w}.
\]
The invariant homology in positive degree $j$ vanishes unless
$j$ is a multiple of $2(d-1)$.
\end{lemma}

\begin{proof}
The boundary datum and its image under the connecting map both
transform by $\epsilon_r$. Duality caps with $[M_r,\partial M_r]$,
on which $s$ also acts by $\epsilon_r$. The two signs cancel.
To compute the invariant part, recall that $H^*(M_r;\R)$ is
generated by the classes $g_{ij}$ of degree $d-1$, pulled back
from $S^{d-1}$ by the pair-direction maps
\cite[Section~8, equation~(101)]{LV}. The reflection acts by $-1$
on the sphere's top cohomology, hence on each $g_{ij}$. It acts
by $(-1)^m$ on degree $m(d-1)$, and the other positive degrees
vanish. Over $\R$ the dual homology representation has the
same signs.
\end{proof}

\subsection{Extension through weight \texorpdfstring{$2d-1$}{2d-1}}
\begin{theorem}\label{thm:weight-range}
Let $d\geq2$, and let
$V=\bigoplus_{a=0}^W V_{(a)}$ be as above,
with $W\leq 2d-1$. Every
$\mathsf L_{\infty,d}$-algebra structure on
$V$ whose primary operations have weight one
extends to an $\En$-algebra structure
recovering those operations along $\iota_d$.
The extension can be chosen so that the
degree-zero binary product vanishes.
\end{theorem}

\begin{proof}
If $W=0$, all primary operations vanish, and we take the zero
action in arities at least two. Otherwise begin with
\eqref{eq:primary-weight-form}, choosing the reflected average above.
If $W\geq2$ we next solve the weight-two problem.
At weight two, the prescribed form has
degree $D_r-1$, and on each face it is
\[
\gamma_i^*K_{r,2}
=\bigl(\operatorname{pr}_1^*\eta_k
       \wedge\operatorname{pr}_2^*\eta_l\bigr)
       (\tau_k\circ_i\tau_l).
\]
The signed fundamental-boundary identity
\eqref{eq:fundamental-boundary} and
$\int_{M_s}\eta_s=1$ give
\[
\int_{\partial M_r}K_{r,2}
=\sum_T\varepsilon(T)
       \tau_{k(T)}\circ_T\tau_{l(T)}
=0.
\]
The last equality is the image of
$d\mu_r$ under the prescribed
$\mathsf L_{\infty,d}$-action.
Its target differential is zero.
Equation~\eqref{eq:top-boundary-integral}
therefore supplies a closed extension.
Average it over permutations and with the reflection character to
obtain $\alpha_{r,2}$ with both prescribed symmetries.

For weights $3\leq w\leq W$, use
Proposition~\ref{prop:weight-obstruction} and
Lemma~\ref{lem:reflection-obstruction}. Whenever the form degree
is nonnegative, the relevant invariant homology vanishes since
$0<w-2<2(d-1)$. Thus at every such stage we may take a closed
extension and impose both symmetries by averaging. In negative
form degree the prescribed data are zero.
Whenever the boundary data are zero,
choose the zero extension, in particular
in arity two and weights at least two.

For a homogeneous $p$-chain $c$, define
\begin{equation}\label{eq:filtered-chain-action}
\rho_r(c)=
\begin{cases}
\displaystyle\int_c \alpha_{r,w},
&w=d(r-1)-p\in\{1,\ldots,W\},\\[4pt]
0,&\text{otherwise}.
\end{cases}
\end{equation}
The unary point acts as the identity.
Each integral has endomorphism degree $-p$.
Closedness and Stokes' theorem prove the
chain-map identity, and the equivariance
of the forms proves equivariance of $\rho$.

For completeness, let $c$ and $c'$ have
dimensions $p$ and $q$, in arities $k,l\geq2$.
Put
\[
a=d(k-1)-p,\qquad b=d(l-1)-q.
\]
Geometric chains have dimensions at most
$D_k,D_l$, so $a,b\geq1$.
If $a+b\leq W$, restriction by
\eqref{eq:weight-boundary-equation} and
integration over the product give
\[
\rho_r(c\circ_i c')
=\rho_k(c)\circ_i\rho_l(c').
\]
Only the summand with factor degrees $p,q$
can contribute to that integral.
If $a+b>W$, the left-hand side is zero
by definition, and the right-hand side
is zero because its weight exceeds $W$.
This also covers the case that one factor
acts by zero because its weight exceeds
$W$. Unary substitutions are identities.
Thus all operadic compositions hold.

On the fundamental chain only weight one
can contribute, and
\[
\rho_r([M_r])=\tau_r.
\]
Finally a zero-chain in $M_2$ has weight
$d$. For $W<d$ it acts by zero automatically.
For $W\geq d$ it acts by $\alpha_{2,d}=0$.
\end{proof}

\subsection{Consequences for secondary polytopes}
\begin{corollary}\label{cor:small-enhancement}
Let $A\subset\R^d$, $d\geq2$, be finite,
affinely spanning, and in affine general
position. If $|A|\leq d+3$, then
$B_A=\gA[1-d]$ admits a nonunital
$\En$-structure recovering its shifted
secondary-polytope operations.
Its degree-zero binary product can be
chosen to vanish.
\end{corollary}

\begin{proof}
Lemma~\ref{lem:small-factorization} verifies
the factorization hypothesis.
The geometric differential is zero,
all generator weights lie in
$\{0,1,2\}$, and the primary operations
raise weight by one by
Lemma~\ref{lem:weight-growth}.
Apply Theorem~\ref{thm:weight-range}.
\end{proof}

\begin{corollary}\label{cor:dimension-range}
Suppose $A\subset\R^d$ satisfies the same
affine general-position assumptions and
also the factorization hypothesis of
Definition~\ref{def:factorization}.
If $|A|\leq3d$, the same enhancement exists.
\end{corollary}

\begin{proof}
The maximum weight is $|A|-d-1\leq2d-1$. Apply
Lemma~\ref{lem:weight-growth} and Theorem~\ref{thm:weight-range}.
Together with Corollary~\ref{cor:small-enhancement}, this proves
Theorem~\ref{thm:intro-geometric}.
\end{proof}

\subsection{Coefficients generated by primary operations}
The extension can be carried out in smaller coefficient spaces.
Let $T_{r,w}\subseteq E_{r,w}$ be the span of evaluations of
rooted trees with $r$ labelled leaves and exactly $w$ internal
vertices. A vertex of valence $k\geq2$ is labelled by $\tau_k$.
The identity $\sum_v(k_v-1)=r-1$ gives
\[
 \sum_v |\tau_{k_v}|=w-d(r-1)=-p_{r,w}.
\]
Every such evaluation has weight $w$, and
\begin{equation}\label{eq:tree-coefficient-composition}
 T_{k,a}\circ_iT_{l,b}\subseteq T_{k+l-1,a+b}.
\end{equation}

\begin{proposition}\label{prop:tree-coefficients}
Under the hypotheses of Theorem~\ref{thm:weight-range}, the
integration action can be chosen with every $\alpha_{r,w}$ valued
in $T_{r,w}$. Moreover,
\[
 T_{r,w}=0\quad\text{if }r<w+1.
\]
If $\tau_k=0$ for $k>R$, then $T_{r,w}=0$ also for
$r>1+w(R-1)$. If only $\tau_2$ is nonzero, then
$T_{r,w}=0$ unless $r=w+1$.
\end{proposition}

\begin{proof}
The primary form has coefficient $\tau_r\in T_{r,1}$.
Inductively, \eqref{eq:tree-coefficient-composition} places the
prescribed boundary form in $T_{r,w}$. Apply the relative extension
lemma with this coefficient space. At weight two its integral
vanishes by the same $L_\infty$ relation; in the higher weights
of Theorem~\ref{thm:weight-range}, the required scalar invariant
homology vanishes. Relabelling preserves $T_{r,w}$, and reflection
acts trivially on coefficients, so both averages stay in this space.
If it is zero, take the zero form.

Every internal vertex has $k_v-1\geq1$. Thus $r-1\geq w$.
If $k_v\leq R$ at every nonzero vertex, then
$r-1\leq w(R-1)$. For binary trees both inequalities are equalities.
\end{proof}

\begin{corollary}\label{cor:polygon-chain-degrees}
Let $A$ be the vertices of a strictly convex $m$-gon, with
$3\leq m\leq6$. Its geometric $E_2$-enhancement can be chosen
so that its only possibly nonzero nonunary components are
\[
 C^{\mathrm{SA}}_{r-1}(\FM_2(r);\R)
       \longrightarrow\End^{1-r}_{B_A}(r),
 \qquad 2\leq r\leq m-2.
\]
\end{corollary}

\begin{proof}
Factorization follows from Lemma~\ref{lem:polygon-factorization}.
For a convex polygon every coarse proper subdivision is a split
along a single diagonal: its cell dual graph is a tree, and any
dissection with more than one diagonal admits a proper coarsening.
Such splits are regular, for instance by a convex hinge height
function along the diagonal. Thus only $\tau_2$ is nonzero.
Proposition~\ref{prop:tree-coefficients} restricts the nonzero forms
to $w=r-1$, of degree $2(r-1)-w=r-1$. Finally
$W=m-3\leq3=2d-1$ bounds their arities by $r\leq m-2$.
\end{proof}

\subsection{A period criterion for the next extension}
Suppose a partial family of reflected closed forms satisfying
\eqref{eq:weight-boundary-equation} has been chosen through weight
$2d-1$. This hypothesis makes sense for any finite maximum weight
$W\geq2d$, before the construction gives a full action.
Write $K_{r,2d}$ for the prescribed next boundary form. It has degree
\[
 p_{r,2d}=d(r-3).
\]
For a pair of labels $i,j$, let $\omega_{ij}$ be the pullback of a
normalized closed volume form on $S^{d-1}$ by the pair-direction
map $M_r\to S^{d-1}$. Products of two such classes span
$H^{2d-2}(M_r;\R)$ \cite[Section~8]{LV}.

\begin{proposition}\label{prop:period-criterion}
For $p_{r,2d}\geq0$, the obstruction to a closed extension of
$K_{r,2d}$ vanishes if and only if
\begin{equation}\label{eq:obstruction-periods}
 \int_{\partial M_r} K_{r,2d}\wedge\omega_{ij}\wedge\omega_{kl}=0
\end{equation}
for every choice of pairs. The integrals take values in the
coefficient space. It suffices to test products representing a
basis of $H^{2d-2}(M_r;\R)$.
\end{proposition}

\begin{proof}
Choose an extension $\widetilde K$ of the boundary form. The
obstruction is $[d\widetilde K]$ in relative degree $p_{r,2d}+1$.
Poincare--Lefschetz duality gives a perfect pairing with absolute
cohomology in degree
$D_r-p_{r,2d}-1=2d-2$. For a closed test form $\zeta$ of that
degree, Stokes' theorem gives
\[
 \int_{M_r}d\widetilde K\wedge\zeta
 =\int_{\partial M_r}K_{r,2d}\wedge\zeta.
\]
Apply this coefficientwise and use the spanning products above.
Vanishing of every pairing is equivalent to vanishing of the
relative class. The closed extension then follows from
Lemma~\ref{lem:relative-form-extension} and can be averaged to
preserve the symmetries.
\end{proof}

With coefficients in $T_{r,2d}$, only arities $r\geq2d+1$ can
contribute, by Proposition~\ref{prop:tree-coefficients}. The next
section compares extension choices within the established range;
Section~\ref{sec:further-directions} describes the first calculation
beyond it.

\section{Homotopies of geometric enhancements}
\label{sec:geometric-homotopies}
We compare the integration actions by constructing their homotopies
on the same configuration spaces. The reflection used for existence
also controls the relative cylinder problem. In this section the
primary operations $\tau_r$ are fixed, while their representing
forms may vary.

\subsection{The relative cylinder complex}
Let $X$ be a compact smooth manifold with corners whose boundary
hypersurfaces are embedded, and put
$R_I=\R[t,dt]$, with $|t|=0$, $|dt|=1$, $d(t)=dt$.
Let $R_{I,\partial I}$ be the kernel of endpoint evaluation
$R_I\to\R\oplus\R$. In degree zero it consists of polynomials
vanishing at both endpoints; in degree one it consists of all
polynomial multiples of $dt$. Integration gives
\[
 H^j(R_{I,\partial I})=
 \begin{cases}\R,&j=1,\\0,&j\neq1.\end{cases}
\]
For degree one this follows because a polynomial $f(t)dt$ has a
primitive vanishing at both endpoints exactly when
$\int_0^1 f(t)dt=0$.

\begin{lemma}\label{lem:polynomial-cylinder}
Suppose closed total $p$-forms are prescribed on
$\partial X\times I$ and on $X\times\{0,1\}$, polynomial in the
interval variable, and their restrictions agree at the intersections.
The obstruction to a closed total extension belongs to
\[
 H^p(X,\partial X;\R)\otimes E,
\]
where $E$ is the finite-dimensional coefficient space.
If this class is zero, a polynomial closed extension exists.
Compatible data for a finite group admit an equivariant extension
whenever they admit an extension.
\end{lemma}

\begin{proof}
The kernel of the combined boundary and endpoint restriction is
\[
 \Omega^*(X,\partial X)\otimes R_{I,\partial I}.
\]
The restriction onto compatible data is surjective. First extend
the endpoint forms by linear interpolation in $t$. Subtract their
boundary values from the prescribed boundary forms; the difference
has zero endpoint restrictions. The coefficientwise extension of
Lemma~\ref{lem:relative-form-extension}, tensored with
$R_{I,\partial I}$, extends this difference to $X$.
The connecting homomorphism therefore places the obstruction in
degree $p+1$ of this kernel. The relative interval calculation and
the Kunneth formula over $\R$ identify this group with
$H^p(X,\partial X)\otimes E$. Subtracting a relative primitive
gives a closed extension when the class vanishes. Finite-group
averaging preserves the prescribed data and polynomial dependence.
\end{proof}

\subsection{Comparison through weight \texorpdfstring{$2d-2$}{2d-2}}
Use the spaces $V$, $E_{r,w}$ and degrees $p_{r,w}$ of
Section~\ref{sec:weight-extension}. Fix a hyperplane reflection $s$.
Let $\alpha^0_{r,w}$ and $\alpha^1_{r,w}$ be two families of closed
forms satisfying the permutation, insertion and reflection equations.
Their primary terms are
\[
 \alpha^j_{r,1}=\eta^j_r\tau_r,\qquad
 \int_{M_r}\eta^j_r=1,\qquad j=0,1,
\]
with $\eta^j_r$ supported in the interior. The existence construction
has zero higher-weight binary forms. The comparison also allows
binary forms satisfying the same closedness and symmetry conditions.

\begin{theorem}\label{thm:reflected-comparison}
Suppose $W\leq2d-2$. The two integration actions are connected by
a polynomial dg-operad homotopy
\[
 \mathcal H:\En\longrightarrow\End_V\otimes R_I
\]
which is constant on the prescribed shifted Lie action. If both
endpoint families have zero higher-weight binary forms, the homotopy
can be chosen with those terms zero.
\end{theorem}

\begin{proof}
For $W=0$ the nonunary actions are zero, so the constant homotopy
suffices. Assume $W\geq1$. We construct total closed forms
\[
 \widetilde\alpha_{r,w}
 \in\bigl(\Omega^*(M_r)\otimes R_I\bigr)^{p_{r,w}}
                      \otimes E_{r,w}
\]
with the specified endpoints. The boundary equations are the
equations of Section~\ref{sec:weight-extension}, applied to total
forms, with the same interval variable on both insertion factors.

For weight one, the compactly supported top form
$\eta^1_r-\eta^0_r$ has integral zero. Compactly supported de Rham
cohomology of the connected oriented interior gives a compactly
supported $(D_r-1)$-form $\beta_r$ such that
\[
 d_M\beta_r=\eta^1_r-\eta^0_r.
\]
Average it with the permutation and reflection characters. With
$dt$ written on the left, set
\begin{equation}\label{eq:primary-cylinder}
 \widetilde\alpha_{r,1}
  =\bigl((1-t)\eta^0_r+t\eta^1_r+dt\wedge\beta_r\bigr)\tau_r.
\end{equation}
This total form is closed. It vanishes tangentially on every
insertion face, and its integral over each top-dimensional slice
equals $\tau_r$.

Suppose the lower weights have been constructed. At weight $w\geq2$,
their compositions give the prescribed total forms on each boundary
face. Associativity and the sign cancellation for disjoint insertions
in Proposition~\ref{prop:weight-obstruction} show that these forms
agree on corners. They are closed and have the endpoint values
prescribed by $\alpha^0_{r,w}$ and $\alpha^1_{r,w}$.
Lemma~\ref{lem:polynomial-cylinder} places the obstruction in
\[
 H^{p_{r,w}}(M_r,\partial M_r;\R)\otimes E_{r,w}.
\]
The data transform by $\epsilon_r=(-1)^{r-1}$ under $s$.
Duality caps with the relative fundamental class, which has the
same character. Thus the obstruction belongs, under this duality, to
\begin{equation}\label{eq:cylinder-obstruction}
 H_{w-1}(M_r;\R)^{s=+1}\otimes E_{r,w}.
\end{equation}
For $2\leq w\leq2d-2$ we have
$1\leq w-1<2(d-1)$, so this invariant homology is zero by
Lemma~\ref{lem:reflection-obstruction}. We obtain a closed extension
and average over both symmetries. Negative total form degrees
carry zero data. The same cylinder argument applies in arity two,
where the insertion boundary is empty. If its higher-weight endpoint
forms are both zero, choose the zero homotopy in those components.

It remains to pass from total forms to an operadic homotopy with
the stated signs. Write the contribution of weight $w$ as
\[
 \widetilde\alpha_{r,w}=a_{r,w}(t)+dt\wedge b_{r,w}(t),
 \qquad p=p_{r,w}.
\]
Its coefficient endomorphisms have degree $-p$; the spatial form
degrees of $a$ and $b$ are $p$ and $p-1$. Closedness is equivalent to
\[
 d_Ma_{r,w}=0,\qquad \partial_t a_{r,w}=d_Mb_{r,w}.
\]
With endomorphisms to the left of interval coefficients, define
the weight-$w$ contribution on a homogeneous $q$-chain by
\begin{equation}\label{eq:integrated-cylinder}
 \mathcal H_{r,w}(c)=
 \begin{cases}
   \displaystyle\int_c a_{r,w}(t),&q=p,\\[4pt]
   \displaystyle(-1)^p\left(\int_c b_{r,w}(t)\right)dt,&q=p-1,\\[4pt]
   0,&\text{otherwise}.
 \end{cases}
\end{equation}
Sum over the weights and let the unary point act as the identity.
The differential on $\End_V\otimes R_I$ is
$d(f\otimes u)=(-1)^{|f|}f\otimes du$. For $q=p$,
Stokes gives
\[
 d\mathcal H_{r,w}(c)
 =(-1)^p\left(\int_c\partial_t a_{r,w}\right)dt
 =(-1)^p\left(\int_{\partial c}b_{r,w}\right)dt
 =\mathcal H_{r,w}(\partial c).
\]
For $q=p+1$ this identity is $\int_{\partial c}a_{r,w}=0$;
in the other degrees it follows from $dt^2=0$ or dimension.

The total boundary equations imply composition. To check the interval
sign, consider outer and inner total forms of degrees $P,Q$.
An outer $dt$ contribution in their product has the factor
$(-1)^{P+Q}$ in \eqref{eq:integrated-cylinder}. On the endomorphism
side, the outer integration gives $(-1)^P$, and moving its $dt$
past the inner degree-$-Q$ endomorphism gives $(-1)^Q$.
An inner $dt$ contribution has an additional wedge sign $(-1)^P$
on the form side, hence total sign $(-1)^Q$, exactly its inner
integration sign. Two $dt$ contributions vanish on both sides.
The terms without $dt$ satisfy the ordinary integration composition
identity. These checks cover the sums over weights; terms of total
weight greater than $W$ vanish. Equivariance follows from the
averaged forms.

Endpoint evaluation gives the prescribed actions. On $[M_r]$ only
the ordinary top-degree primary term can contribute. Its integral
is $\tau_r$ for every $t$, by \eqref{eq:primary-cylinder}. Thus
$\mathcal H\iota_d=\tau\otimes1$.
\end{proof}

\begin{corollary}\label{cor:configuration-comparison}
Let $A\subset\R^d$ satisfy the affine and factorization hypotheses,
with $|A|\leq3d-1$. The reflection-compatible actions supplied by
Theorem~\ref{thm:weight-range}, including their choices of primary
forms, are polynomially homotopic relative to the fixed shifted Lie
structure. This proves Proposition~\ref{prop:intro-comparison}.
\end{corollary}

\begin{proof}
The maximum weight is $|A|-d-1\leq2d-2$. Apply
Theorem~\ref{thm:reflected-comparison} to $V=B_A$.
\end{proof}

\subsection{An explicit comparison in weight two}
The low-weight case also admits a direct formula with the primary
forms fixed. We record it independently of the reflection hypothesis.
\label{sec:weight-two-homotopy}

Let $d\geq2$ and
\[
V=V_{(0)}\oplus V_{(1)}\oplus V_{(2)}
\]
be finite-dimensional, with zero differential.
Fix a $\mathsf L_{\infty,d}$-action whose
primary operations $\tau_r$ have weight one,
and fix the normalized forms $\eta_r$.
Write
\[
q_r=D_r-1=d(r-1)-2,\qquad
E_r=E_{r,2}=\End^{-q_r}_{(2)}(V)(r).
\]
Let $\Omega_r^0$ and $\Omega_r^1$ be two
families of closed equivariant $E_r$-valued
$q_r$-forms satisfying the same boundary
equation
\begin{equation}\label{eq:fixed-weight-two-boundary}
\gamma_i^*\Omega_r^j
=\bigl(\operatorname{pr}_1^*\eta_k
       \wedge\operatorname{pr}_2^*\eta_l\bigr)
       (\tau_k\circ_i\tau_l),\qquad j=0,1.
\end{equation}
As in Section~\ref{sec:weight-extension},
assume $\Omega_2^0=\Omega_2^1=0$.
These families determine actions
$\rho^0,\rho^1$ by
\eqref{eq:filtered-chain-action}.

The differences
\[
\Delta_r=\Omega_r^1-\Omega_r^0
\]
are closed relative forms: their tangential
restrictions to all boundary faces vanish.
They therefore determine classes
\begin{equation}\label{eq:relative-difference-class}
[\Delta_r]\in
H^{q_r}(M_r,\partial M_r;\R)\otimes E_r.
\end{equation}
These classes measure relative exactness of
the chosen differences. We will use their
vanishing as a sufficient condition for
an operadic homotopy.

\begin{lemma}\label{lem:relative-difference-primitive}
If $[\Delta_r]=0$, there is an equivariant
relative form
\[
\beta_r\in\Omega^{q_r-1}(M_r,\partial M_r;E_r),
\qquad d\beta_r=\Delta_r.
\]
For $q_r=0$, this assertion means that
$\Delta_r=0$ and $\beta_r=0$.
\end{lemma}

\begin{proof}
This is the definition of a zero class in
the relative de Rham complex used in
Lemma~\ref{lem:relative-form-extension}.
Average a primitive over $S_r$, using
pullback on forms and the inverse coefficient
action. Since $\Delta_r$ is equivariant,
the average is still a primitive.
Its boundary restriction remains zero.
In degree zero there are no boundaries
from negative-degree forms, so a zero
cohomology class is the zero form.
\end{proof}

\subsubsection{The explicit homotopy}
Put $R=\R[t,dt]$, with $|t|=0$,
$|dt|=1$, $d(t)=dt$, and $(dt)^2=0$.
Endomorphisms are written on the left
of coefficients in $R$.
The differential and composition conventions are
\[
\begin{aligned}
d(f\otimes a)&=(-1)^{|f|}f\otimes da,\\
(f\otimes a)\circ_i(g\otimes b)
  &=(-1)^{|a||g|}(f\circ_i g)\otimes ab.
\end{aligned}
\]
The first identity uses $d_V=0$.
By a polynomial homotopy we mean a
dg-operad map to $\End_V\otimes R$
with specified evaluations at $t=0,1$.

\begin{proposition}\label{prop:weight-two-homotopy}
Under the fixed-data assumptions above,
suppose $[\Delta_r]=0$ for every $r$.
Then there is a polynomial homotopy
\[
\mathcal H:\En\longrightarrow\End_V\otimes R
\]
with endpoints $\rho^0,\rho^1$, such that
\[
(\mathcal H\circ\iota_d)(\mu_r)=\tau_r\otimes1.
\]
In particular the prescribed
$\mathsf L_{\infty,d}$-action remains fixed.
\end{proposition}

\begin{proof}
Choose the primitives of
Lemma~\ref{lem:relative-difference-primitive},
taking $\beta_2=0$.
Set $\Omega_r(t)=(1-t)\Omega_r^0+t\Omega_r^1$.
For a homogeneous $p$-chain $c$, define
\begin{equation}\label{eq:weight-two-homotopy}
\mathcal H_r(c)=
\begin{cases}
\displaystyle\left(\int_c\eta_r\right)
       \tau_r\otimes1,
     &p=D_r,\\[5pt]
\displaystyle\int_c\Omega_r(t),
     &p=q_r,\\[5pt]
\displaystyle(-1)^{q_r}
       \left(\int_c\beta_r\right)\otimes dt,
     &p=q_r-1,\\[5pt]
0,&\text{otherwise}.
\end{cases}
\end{equation}
The second line is interpreted coefficientwise
in $E_r\otimes\R[t]$. A line with a negative
chain dimension is omitted.
The unary point acts by
$\operatorname{id}_V\otimes1$.

The first two nonzero lines have endomorphism
degree $-p$. In the third line the
endomorphism degree is $-q_r$, and $dt$
adds one, giving total degree $-p$.
For $p=q_r$, Stokes' theorem gives
\[
\begin{aligned}
d\,\mathcal H_r(c)
&=(-1)^{q_r}\left(\int_c\Delta_r\right)\otimes dt\\
&=(-1)^{q_r}
     \left(\int_{\partial c}\beta_r\right)\otimes dt
 =\mathcal H_r(\partial c).
\end{aligned}
\]
For $p=D_r=q_r+1$, the left-hand side
is zero because the primary form is fixed.
The right-hand side is
$\int_{\partial c}\Omega_r(t)=0$ by closedness.
For $p=q_r-1$, the differential of the
$dt$ term is zero and the boundary acts
by zero. Lower degrees vanish on both
sides. If $q_r=0$, the difference is zero
and the second line is constant in $t$,
so the same conclusion holds.
Equivariance follows from that of all
the forms.

It remains to check compositions.
Take chains $c,c'$ in arities $k,l\geq2$
and put $r=k+l-1$. Their total dimension
is at most
$D_k+D_l=D_r-1=q_r$.
If both have top dimension, the left-hand
side of the composition identity is
the integral of $\Omega_r(t)$ on the
insertion face. Equation
\eqref{eq:fixed-weight-two-boundary}
is independent of $t$ and gives exactly
the composition of the primary operations.

If their dimensions sum to $q_r-1$,
the left-hand side is the $dt$ term in
\eqref{eq:weight-two-homotopy}.
It vanishes because $\beta_r$ restricts
to zero on every face.
For a smaller sum it vanishes by definition.
On the right-hand side, unless both
chains have top dimension, any nonzero
factor of lower dimension has weight two.
The other nonunary factor has weight
at least one. Their composite has
weight at least three and is zero on $V$.
This remains true for the $dt$ terms,
regardless of the tensor-product sign.
Unary substitutions satisfy the identity
axiom. These cases prove all operadic
composition identities.

Evaluation at either endpoint kills $dt$
and recovers $\rho^0$ or $\rho^1$.
Finally, the normalized primary form gives
$\mathcal H_r([M_r])=\tau_r\otimes1$.
Since the $\mu_r$ generate the underlying
graded operad, the homotopy is fixed on
the entire image of $\iota_d$.
\end{proof}

\subsubsection{Vanishing of the relative difference classes}
\begin{theorem}\label{thm:weight-two-choice-independence}
If $d\geq3$, any two weight-two extensions
with the fixed primary forms specified
above are polynomially homotopic relative
to the prescribed
$\mathsf L_{\infty,d}$-action.
\end{theorem}

\begin{proof}
Poincare--Lefschetz duality gives
\[
H^{q_r}(M_r,\partial M_r;\R)
\cong H_{D_r-q_r}(M_r;\R)
=H_1(M_r;\R).
\]
This group is zero for $d\geq3$, by the
configuration-space cohomology calculation
used in Lemma~\ref{lem:reflection-obstruction}.
Thus all classes
\eqref{eq:relative-difference-class}
vanish. Apply
Proposition~\ref{prop:weight-two-homotopy}.
\end{proof}

\begin{corollary}\label{cor:small-fixed-form-homotopy}
Let $d\geq3$, and let $A\subset\R^d$ be
affinely spanning and in affine general
position, with $|A|\leq d+3$.
Fix the normalized primary forms.
The actions on $B_A$ constructed in
Corollary~\ref{cor:small-enhancement}
are independent, up to the polynomial
homotopy above, of the choices of
weight-two correction forms.
\end{corollary}

\begin{proof}
The weights of $B_A$ lie in
$\{0,1,2\}$. Apply
Theorem~\ref{thm:weight-two-choice-independence}.
\end{proof}

\begin{corollary}
For $d=2$, with fixed reflection-compatible primary forms, any
two weight-two corrections satisfying \eqref{eq:reflection-type}
give polynomially homotopic actions on a space of weights
$0,1,2$, relative to the prescribed shifted Lie operations.
\end{corollary}

\begin{proof}
The relative difference class transforms by $\epsilon_r$.
Capping with the relative fundamental class identifies it with
an element of $H_1(M_r;\R)^{s=+1}\otimes E_r$, by the same
duality calculation as in Lemma~\ref{lem:reflection-obstruction}.
This invariant space vanishes. Apply
Proposition~\ref{prop:weight-two-homotopy}.
\end{proof}

\section{Examples and the first nonzero composites}
\label{sec:examples}

We give three circuit examples, a complete operation table for a
convex pentagon, and the chain operations required by a convex hexagon. For a simplex with vertices
labelled by $i_1,\ldots,i_s$, write
$t_{i_1\cdots i_s}$ for its generator in $B_A$.
The generator associated to a whole circuit
is denoted by $q$.
We use the sign normalization at the end of
Section~\ref{sec:secondary-polytopes}.
Displayed values are extended using \eqref{eq:tau-symmetry}.
For the circuit examples and the complete pentagon table,
undisplayed basis tuples have value zero. In the hexagon example,
the primary maps are all the subdivision operations of
Section~\ref{sec:secondary-polytopes}; the displayed fan illustrates
one of their nonzero composites.

\subsection{A convex quadrilateral}
Let
\[
a_1=(0,0),\quad a_2=(1,0),\quad
a_3=(1,1),\quad a_4=(0,1).
\]
The affine relation is
$-a_1+a_2-a_3+a_4=0$.
The two triangulations are
\[
\mathcal T_+=\{123,134\},\qquad
\mathcal T_-=\{124,234\},
\]
corresponding to the two diagonals.
Using Euclidean area, their GKZ vectors are
\[
\phi_{\mathcal T_+}=(1,\tfrac12,1,\tfrac12),
\qquad
\phi_{\mathcal T_-}=(\tfrac12,1,\tfrac12,1).
\]
All four triangles are empty. Hence
\[
B_A^2=\R\{t_{123},t_{124},t_{134},t_{234}\},
\qquad B_A^3=\R q.
\]
Orient the secondary interval from
$\mathcal T_-$ to $\mathcal T_+$.
The only nonzero primary operation is
the degree-$-1$ binary map
\[
\tau_2(t_{123},t_{134})=q,\qquad
\tau_2(t_{124},t_{234})=-q.
\]
These values are symmetric under exchanging
the two degree-two inputs. Every operation
with an input $q$ vanishes.

The full $E_2$-action is
\[
\rho_2(c)=\left(\int_c\eta_2\right)\tau_2
\quad\text{for }c\in C^{\mathrm{SA}}_1(S^1;\R),
\]
and it is zero in all other nonunary
arities and chain degrees. The unary
point acts as the identity.
The fundamental circle acts by $\tau_2$;
zero-chains act by zero. Thus the
binary operation in this example is
the degree-$-1$ operation, while the
degree-zero product vanishes.

\subsection{A triangle with an interior point}
Take
\[
a_1=(0,0),\quad a_2=(3,0),\quad
a_3=(0,3),\quad a_4=(1,1).
\]
The relation $a_1+a_2+a_3-3a_4=0$ gives
\[
\mathcal T_\star=\{124,134,234\},
\qquad
\mathcal T_0=\{123\}.
\]
Their GKZ vectors are
\[
\phi_{\mathcal T_\star}=(3,3,3,\tfrac92),
\qquad
\phi_{\mathcal T_0}=(\tfrac92,\tfrac92,\tfrac92,0).
\]
The outer triangle with only its three
vertices marked is not geometric:
it contains $a_4$. The geometric
graded space is therefore
\[
B_A^2=\R\{t_{124},t_{134},t_{234}\},
\qquad B_A^3=\R q.
\]
Orient the secondary interval from
$\mathcal T_0$ to $\mathcal T_\star$.
Its full cellular differential is
\[
d(v_A)=v_{124}v_{134}v_{234}-v_{123}.
\]
After passing to the geometric algebra,
the unary term has no geometric input.
The only nonzero primary map is
\begin{equation}\label{eq:interior-ternary}
\tau_3(t_{124},t_{134},t_{234})=q,
\qquad |\tau_3|=-3.
\end{equation}
In particular, $\tau_2=0$.

Theorem~\ref{thm:central-extension} gives
\[
\rho_3(c)=\left(\int_c\eta_3\right)\tau_3
\quad\text{for }c\in
C^{\mathrm{SA}}_3(\FM_2(3);\R),
\]
with every other nonunary chain degree
and arity acting by zero. The degree
check is $2+2+2-3=3$.
The fundamental chain of $\FM_2(3)$
has boundary, so \eqref{eq:interior-ternary}
is a chain operation rather than an
operation represented by its absolute
fundamental homology class.

The induced binary product and binary
bracket on cohomology are both zero.
Nevertheless, the ternary operation
cannot be eliminated by a change of
minimal $L_\infty$ coordinates.

\begin{proposition}\label{prop:essential-ternary}
Let $\mathfrak g_A=B_A[1]$ be the
$L_\infty$-algebra of this example.
If a minimal $L_\infty$-algebra
$\mathfrak h$ is related to
$\mathfrak g_A$ by an $L_\infty$
quasi-isomorphism in either direction,
then its binary bracket is zero and
its ternary bracket is nonzero.
In particular, $\mathfrak g_A$ is not
$L_\infty$-isomorphic to an abelian
graded algebra.
\end{proposition}

\begin{proof}
Use the symmetric coalgebra convention
of Section~\ref{sec:secondary-polytopes}.
The coderivation $\mathcal Q$ of
$\mathfrak g_A[1]$ has
$\mathcal Q_1=\mathcal Q_2=0$ and
$\mathcal Q_3\neq0$.
Suppose first that
$F:\mathfrak g_A\to\mathfrak h$
is an $L_\infty$ quasi-isomorphism,
and denote the target coderivation
by $\mathcal Q'$.
Both algebras are minimal, so $F_1$
is an isomorphism of graded spaces.
The length-two part of
$F\mathcal Q=\mathcal Q'F$ gives
$\mathcal Q'_2F_1^{\otimes2}=0$,
hence $\mathcal Q'_2=0$.
The length-three part then reduces to
\[
F_1\mathcal Q_3
=\mathcal Q'_3F_1^{\otimes3}.
\]
Its left-hand side is nonzero, so
$\mathcal Q'_3\neq0$. The same equations
with source and target interchanged
prove the assertion for a morphism
in the other direction.
\end{proof}

Thus the first nonzero higher operation
contains information that is absent
from the binary cohomology operations.
The proposition concerns the specified
minimal $L_\infty$ structures and uses
no classification of all $E_2$-algebras.

\subsection{A triangular bipyramid}
Let $A\subset\R^3$ consist of
\[
\begin{aligned}
a_1&=(0,0,0),&a_2&=(4,0,0),&a_3&=(0,4,0),\\
a_4&=(1,1,2),&a_5&=(1,1,-2).
\end{aligned}
\]
Its affine relation is
$2a_4+2a_5=2a_1+a_2+a_3$.
The two triangulations are
\[
\mathcal T_2=\{1234,1235\},\qquad
\mathcal T_3=\{1245,1345,2345\}.
\]
The first glues two tetrahedra along
the triangle $123$; the second glues
three tetrahedra around the edge $45$.
Their total volumes are both $32/3$,
and their GKZ vectors are
\[
\begin{aligned}
\phi_{\mathcal T_2}
  &=(\tfrac{32}3,\tfrac{32}3,\tfrac{32}3,
     \tfrac{16}3,\tfrac{16}3),\\
\phi_{\mathcal T_3}
  &=(\tfrac{16}3,8,8,\tfrac{32}3,\tfrac{32}3).
\end{aligned}
\]
All five tetrahedra are empty.
The underlying space is six-dimensional:
\[
\begin{aligned}
B_A^3&=\R\{t_{1234},t_{1235},t_{1245},
                  t_{1345},t_{2345}\},\\
B_A^4&=\R q.
\end{aligned}
\]
Orient the interval with boundary
$\mathcal T_2-\mathcal T_3$.
For $d=3$, the suspension signs
$\epsilon_{3,2}$ and $\epsilon_{3,3}$
are both $-1$. Taking $q$ to be the
negative of the canonical shifted
dual generator gives
\begin{equation}\label{eq:bipyramid-operations}
\tau_2(t_{1234},t_{1235})=q,\qquad
\tau_3(t_{1245},t_{1345},t_{2345})=-q.
\end{equation}
Their degrees are respectively $-2$
and $-5$, as checked by
\[
3+3-2=4,\qquad 3+3+3-5=4.
\]
Both displayed values are unchanged
by permuting the inputs. For odd $d$,
this is compatible with the operadic
sign character: permuting degree-three
inputs contributes the same sign as
the permutation character of the
configuration-space orientation.

The full action has two nonzero
nonunary components:
\[
\begin{aligned}
\rho_2(c)&=\left(\int_c\eta_2\right)\tau_2,
   &&\dim c=2,\quad M_2=S^2,\\
\rho_3(c)&=\left(\int_c\eta_3\right)\tau_3,
   &&\dim c=5,\quad M_3=\FM_3(3).
\end{aligned}
\]
All remaining nonunary chain operations
are zero. Here $\eta_2$ is a normalized
sphere volume form and $\eta_3$ is
chosen as in Lemma~\ref{lem:normalized-forms}.
In both arities it transforms by
$\sgn(\sigma)$. Theorem~\ref{thm:central-extension}
proves the composition identities,
and Proposition~\ref{prop:form-homotopy}
controls the choice of forms.

\subsection{The first nonzero composite: a convex pentagon}
Let $A=\{a_1,\ldots,a_5\}\subset\R^2$
be the vertices of a strictly convex
pentagon, labelled cyclically and in
affine general position. One example is
\[
(0,0),\quad(2,0),\quad(3,1),\quad
(1,3),\quad(-1,1).
\]
Every subconfiguration of at least
three points is geometric. Thus
\[
\dim B_A^2=10,\qquad
\dim B_A^3=5,\qquad
\dim B_A^4=1.
\]
Write $q_{ijkl}$ for a quadrilateral
generator and $p=p_{12345}$ for the
pentagon generator.

The secondary polytope of a convex
pentagon is a pentagon. Its five
facets correspond to subdivisions
along a diagonal into a triangle
and a quadrilateral. Together with
the circuit operations on its
quadrilaterals, these are all the
primary incidence operations.
Consequently $\tau_r=0$ for $r\neq2$.

Choose the quadrilateral generators so that, for every
$1\leq i<j<k<l\leq5$,
\begin{equation}\label{eq:pentagon-quadrilateral-table}
 \tau_2(t_{ijk},t_{ikl})=q_{ijkl},\qquad
 \tau_2(t_{ijl},t_{jkl})=-q_{ijkl}.
\end{equation}
Choose the sign of $p$ by
$\tau_2(q_{1234},t_{145})=p$. The remaining coefficients are as follows.

\begin{proposition}[The pentagon operation table]
\label{prop:pentagon-table}
With these choices, the triangle--quadrilateral operations are
\[
\begin{array}{c|c|c}
 \text{triangle}&\text{quadrilateral}&\tau_2\text{ value}\\ \hline
 t_{123}&q_{1345}&-p\\
 t_{234}&q_{1245}& p\\
 t_{345}&q_{1235}&-p\\
 t_{145}&q_{1234}& p\\
 t_{125}&q_{2345}& p
\end{array}
\]
Together with \eqref{eq:pentagon-quadrilateral-table}, symmetry
on these inputs, and zero on all other basis pairs, this is the
complete primary operation. All $\tau_r$ for $r\neq2$ are zero.
\end{proposition}

\begin{proof}
Every facet of the secondary pentagon splits the polygon along
one diagonal, so the five displayed pairs are precisely the
possible triangle--quadrilateral inputs. Each coefficient has
absolute value one by the cellular incidence construction.
Let $s_1,s_2,s_3,s_4,s_5$ denote their signs in the displayed order.
The normalization gives $s_4=1$.

In the symmetric coalgebra convention, the triangle generators have
degree zero after applying $s_2:B_A\to\mathfrak g_A[1]$.
Thus the three-input equation on triangles is the unsigned sum
of the three possible two-vertex evaluations. For the five
triangulations, it gives
\[
\begin{array}{c|c}
 \text{triangles}&\text{equation}\\ \hline
 123,134,145&s_4+s_1=0\\
 123,135,345&s_3-s_1=0\\
 125,235,345&-s_3-s_5=0\\
 125,234,245&s_5-s_2=0\\
 124,145,234&s_2-s_4=0.
\end{array}
\]
Hence $(s_1,s_2,s_3,s_4,s_5)=(-1,1,-1,1,1)$.
These are all possible nonzero two-vertex evaluations: a nonzero
one glues three triangles to the whole pentagon, hence determines
one of its five triangulations. If one of the three inputs already
has positive weight, its output weight would exceed two, so that
evaluation is zero. The formulas therefore give the full table
and explicitly verify its quadratic identity.
\end{proof}

For example,
\begin{equation}\label{eq:pentagon-composite}
 \tau_2\bigl(\tau_2(t_{123},t_{134}),t_{145}\bigr)=p,
 \qquad
 \tau_2\bigl(\tau_2(t_{134},t_{145}),t_{123}\bigr)=-p.
\end{equation}
The third composite on this triple is zero. The two nonzero
terms are the two ways of reaching its triangulation through
facets of the secondary pentagon.

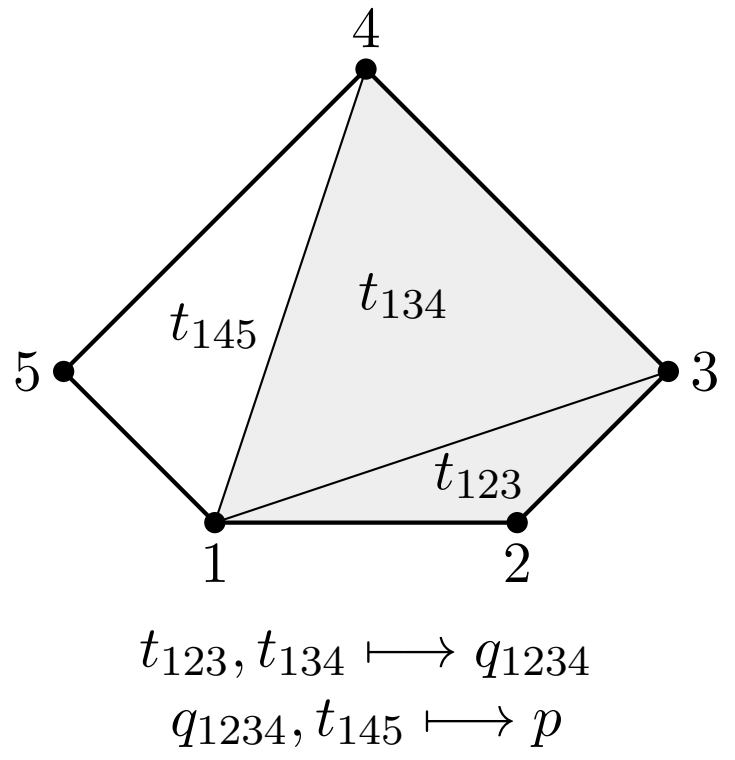
\begin{figure}[htbp]
\centering
\begin{tikzpicture}[scale=1,every node/.style={font=\small}]
 \coordinate (1) at (0,0); \coordinate (2) at (2,0);
 \coordinate (3) at (3,1); \coordinate (4) at (1,3);
 \coordinate (5) at (-1,1);
 \fill[gray!13] (1)--(2)--(3)--(4)--cycle;
 \draw[thick] (1)--(2)--(3)--(4)--(5)--cycle;
 \draw (1)--(3) (1)--(4);
 \foreach \i in {1,2,3,4,5} \fill (\i) circle (2pt);
 \node[below] at (1) {$1$}; \node[below] at (2) {$2$};
 \node[right] at (3) {$3$}; \node[above] at (4) {$4$};
 \node[left] at (5) {$5$};
 \node at (1.75,.31) {$t_{123}$};
 \node at (1.25,1.5) {$t_{134}$};
 \node at (0,1.3) {$t_{145}$};
 \node[align=center] at (1,-1.1)
 {$t_{123},t_{134}\longmapsto q_{1234}$\\
  $q_{1234},t_{145}\longmapsto p$};
\end{tikzpicture}
\caption{A nonzero two-step gluing. The two successive operations
determine a nonzero operation on an insertion face of $\FM_2(3)$.}
\label{fig:pentagon-gluing}
\end{figure}

The following is chain action for the pentagon.
Here $M_2=S^1$ and $\dim M_3=3$. Normalize the circle form
$\eta_2$ so that $\int_{S^1}\eta_2=1$.
The tree-coefficient construction supplies a closed
$T_{3,2}$-valued two-form $\Omega_3$ on $M_3$ with
\begin{equation}\label{eq:pentagon-face-form}
 \gamma_i^*\Omega_3
   =(\operatorname{pr}_1^*\eta_2\wedge
     \operatorname{pr}_2^*\eta_2)(\tau_2\circ_i\tau_2).
\end{equation}
The integral of these prescribed boundary forms is zero by the
quadratic identity verified above. Lemma~\ref{lem:relative-form-extension}
therefore gives the closed extension.

The full nonunary action is
\[
\begin{aligned}
 \rho_2(c)&=\left(\int_c\eta_2\right)\tau_2,
     &&\dim c=1,\\
 \rho_3(c)&=\int_c\Omega_3,
     &&\dim c=2,
\end{aligned}
\]
with all other arities and chain dimensions acting by zero.
The unary point acts as the identity. In particular,
\[
 \rho_3([M_2]\circ_1[M_2])=\tau_2\circ_1\tau_2\neq0,
 \qquad \rho_3([M_3])=0.
\]
Equation~\eqref{eq:pentagon-face-form} checks compositions of
the two binary circle operations. Every composite involving
$\Omega_3$ and a further nonunary operation has weight at least
three and is zero. This realizes the complete chain action in
the first example with a nonzero nested gluing. Choosing the
reflection-compatible forms, Corollary~\ref{cor:configuration-comparison}
compares all the choices in this construction.

\subsection{A convex hexagon}
\label{ex:hexagon}
Let $A$ be six cyclically labelled vertices of a strictly convex
hexagon. Its geometric algebra has
\[
 \dim B_A^2=20,\qquad \dim B_A^3=15,\qquad
 \dim B_A^4=6,\qquad \dim B_A^5=1.
\]
Indeed, every subset of at least three vertices is a geometric
subpolygon, and a subset of size $k$ gives a generator of degree
$k-1$. The maximum weight is three.

By Corollary~\ref{cor:polygon-chain-degrees}, choose the only
possibly nonzero forms to be
\[
 \alpha_{2,1}=\eta_2\tau_2,\qquad
 \alpha_{3,2}=\Omega_3,\qquad
 \alpha_{4,3}=\Theta_4,
\]
of degrees one, two and three. The weight-two forms satisfy
\eqref{eq:pentagon-face-form} with the hexagon coefficient spaces.
The boundary values of $\Theta_4$ are the appropriate insertions
\[
 \alpha_{2,1}\star_i\Omega_3,
 \qquad \Omega_3\star_i\alpha_{2,1},
\]
on the faces $M_2\times M_3$ and $M_3\times M_2$ of $M_4$.
The corner values agree by associativity and the sign calculation
of Proposition~\ref{prop:weight-obstruction}.

Since $\dim M_4=5$, the obstruction to a closed three-form
extension is in $H^4(M_4,\partial M_4)$ with coefficients.
Reflection reduces it to $H_1(M_4)^{s=+1}$, which is zero.
Thus $\Theta_4$ exists. The resulting nonunary action is
\[
\begin{array}{c|c|c}
 \text{arity}&\text{chain dimension}&\text{action}\\ \hline
 2&1&(\int_c\eta_2)\tau_2\\
 3&2&\int_c\Omega_3\\
 4&3&\int_c\Theta_4
\end{array}
\]
and every other component is zero.

Let $h=h_{123456}$ be the hexagon generator. The fan triangulation
has successive gluings
\[
\begin{aligned}
 \tau_2(t_{123},t_{134})&=q_{1234},\\
 \tau_2(q_{1234},t_{145})&=p_{12345},\\
 \tau_2(p_{12345},t_{156})&=h,
\end{aligned}
\]
after choosing the displayed output orientations successively.
Thus the corresponding three-dimensional corner chain
$[M_2]\circ[M_2]\circ[M_2]$ in $M_4$ acts nontrivially.
The fundamental chain of $M_4$ has dimension five and acts by
$\tau_4=0$. The extension supplies the lower-dimensional operation
required by the three-step subdivision.

\section{Further directions}
\label{sec:further-directions}

\subsection{Higher weights}
Reflection restricts the possible obstruction weights to
\[
 w-2\in2(d-1)\mathbb Z_{\geq0}.
\]
Weight two is settled by the shifted Lie identities. The next
weight is $2d$, with obstruction in
$H_{2d-2}(M_r;\R)^{s=+1}$ with endomorphism coefficients.
Evaluating this class, allowing the choices made in the preceding
weight to vary, is the next step in extending the geometric
construction beyond the range of Theorem~\ref{thm:weight-range}.

The convex heptagon is the first planar example. With coefficients
in the tree spaces, the remaining form is a four-form on $M_5$,
whose boundary values are the three types
\[
 \alpha_{2,1}\star\alpha_{4,3},\qquad
 \alpha_{3,2}\star\alpha_{3,2},\qquad
 \alpha_{4,3}\star\alpha_{2,1}.
\]
The period criterion of Proposition~\ref{prop:period-criterion}
gives a finite cohomological calculation for this problem.
The configuration-space Poincare polynomial
\[
 \prod_{j=1}^4(1+jt)=1+10t+35t^2+50t^3+24t^4
\]
shows that $H^2(M_5;\R)$ has dimension $35$ before coefficient
and symmetry reductions \cite[Section~8]{LV}. Determining the
resulting periods and their dependence on $\alpha_{4,3}$ is left
for further work.

\subsection{Associative and mixed structures}
The constructions in this paper concern the geometric
$L_\infty$-algebra. The algebra of the infrared also includes directed
$A_\infty$-algebras and mixed operations relating the Lie and
associative structures \cite{KKS}. Enhancing these structures
compatibly is left for future work. A geometric approach would
require compatible operations with interior and boundary inputs.
In dimension two, the Swiss-cheese operad of \cite{Voronov}
provides a setting for the required compatibility. The
factorization-algebra language of \cite{CG} offers a local
formulation of the same question.
Extension to coefficient systems is another direction for the
method developed here.

\subsection{Variation of the factorization condition}
Corollary~\ref{cor:ambient-factorization-test} reduces the
factorization hypothesis to the regular subdivisions of $A$.
In dimension two, Corollary~\ref{cor:kss-factorization-equivalence}
identifies it with the condition of \cite{KSS}. Studying the
genericity of this condition for $d\geq3$ is a separate problem:
the rank criterion for a fixed subdivision must be combined with
control of which subdivisions remain regular as $A$ varies.

\bibliographystyle{alpha}
\bibliography{references}
\end{document}